\documentclass[a4paper, 11pt]{amsart}

\usepackage[T1]{fontenc}
\usepackage[utf8]{inputenc}
\usepackage[english]{babel}
\usepackage{lmodern}

\usepackage{cite}
\usepackage{soul}

\usepackage[left=2.2cm, right=2.2cm]{geometry}

\usepackage{amsmath, amsthm, amssymb, amsfonts}
\usepackage{mathtools}
\usepackage{mathrsfs}
\usepackage{mathdots}
\usepackage{accents}

\usepackage{enumitem} 
\usepackage{array, booktabs}

\usepackage{graphicx}
\usepackage{tikz}
\usetikzlibrary{matrix,arrows,decorations.pathmorphing}
\usepackage[all]{xy}

\usepackage[usenames,svgnames,dvipsnames]{xcolor}
\usepackage{hyperref}
\hypersetup{
    colorlinks=true,
    urlcolor=Black,
    linkcolor=DarkBlue,
    citecolor=DarkBlue,
}

\numberwithin{equation}{section}

\theoremstyle{plain} 
\newtheorem{thm}{Theorem}[section]

\newtheorem{lemma}[thm]{Lemma}
\newtheorem{cor}[thm]{Corollary}

\newtheorem{mainthm}{Theorem}
\newtheorem{maincor}{Corollary}

\theoremstyle{definition}
\newtheorem{Def}[thm]{Definition}

\newtheorem{ex}[thm]{Example}

\newtheorem{rmk}[thm]{Remark}

\newcommand{\thickbar}[1]{\accentset{\rule{.55em}{1.10pt}}{#1}}
\newcommand{\biq}{\mathbin{/\mkern-6mu/}}          
\renewcommand{\d}[1]{\mathop{}\!\mathrm{d}#1}       

\newcommand{\qandq}{\quad\text{and}\quad}

\newcommand{\Z}{\mathbb{Z}}      
\newcommand{\Q}{\mathbb{Q}}      
\newcommand{\R}{\mathbb{R}}      
\newcommand{\C}{\mathbb{C}}      
\newcommand{\CP}{\mathbb{C}P}    
\newcommand{\HP}{\mathbb{H}P}    

\DeclareMathOperator{\GL}{GL}    
\DeclareMathOperator{\On}{O}
\DeclareMathOperator{\SO}{SO}    
\DeclareMathOperator{\Un}{U}     
\DeclareMathOperator{\SU}{SU}    
\DeclareMathOperator{\Sp}{Sp}    

\DeclareMathOperator{\Ad}{Ad}        
\DeclareMathOperator{\Aut}{Aut}      
\DeclareMathOperator{\Isom}{Isom}    
\DeclareMathOperator{\Span}{span}    
\DeclareMathOperator{\diag}{diag}    
\DeclareMathOperator{\tr}{tr}        
\DeclareMathOperator{\gr}{gr}        
\DeclareMathOperator{\pr}{pr}        
\DeclareMathOperator{\id}{\mathrm{Id}} 
\DeclareMathOperator{\low}{\mathrm{low}} 
\DeclareMathOperator{\irr}{\mathrm{irr}}   
\DeclareMathOperator{\Sym}{\mathrm{Sym}}   
\DeclareMathOperator{\rank}{rank}    
\DeclareMathOperator{\cohom}{cohom}  

\title{Positively curved orbifolds from infinitesimally polar actions on $\mathrm{Sp}(2)$}

\author{Nicolas Caballero}
\address{Instituto de Matem\'atica e Estat\'istica, Universidade de S\~ao Paulo, Rua do Mat\~ao 1010, Cidade Universit\'aria, 05508-090 S\~ao Paulo - Brasil}
\email{nrrcaballero@usp.br}

\author{Daniel Rotmeister}
\address{Instituto de Matem\'atica e Estat\'istica, Universidade de S\~ao Paulo, Rua do Mat\~ao 1010, Cidade Universit\'aria, 05508-090 S\~ao Paulo - Brasil}
\email{daniel.rotmeister@ime.usp.br}

\date{November 28, 2026}
\thanks{}

\begin{document}

\begin{abstract}
     In this paper, we provide a complete classification of the infinitesimally polar actions on the Lie group $\Sp(2)$ equipped with a bi-invariant Riemannian metric, which consequently yields the classification of all its orbifold biquotients. Among these, we explicitly classify all closed connected subgroups that act almost freely on $\Sp(2)$. As a consequence, via the universal double covering $\Sp(2) \to \SO(5)$, we also deduce the complete classification of all infinitesimally polar actions and orbifold biquotients of $\SO(5)$. Furthermore, we construct a new family of $6$-dimensional Riemannian orbifolds and a $4$-dimensional bad Riemannian orbifold homeomorphic to $S^4$, all with positive sectional curvature and a non-trivial singular locus.
\end{abstract} 

\maketitle


\section{Introduction}

A fundamental problem in Riemannian geometry is the study of spaces with positive sectional curvature. Compared to negative curvature, relatively few examples and obstruction theorems are known. The classical examples are given by the compact rank-one symmetric spaces $S^n$, $\CP^n$, $\HP^n$, and $\mathbb{O}P^2$. In the homogeneous setting, simply connected spaces admitting a homogeneous metric with positive curvature have been classified by Berger \cite{berger} in the normal homogeneous case, and by Wallach and B\'erard-Bergery in even and odd dimensions, respectively. In the non-negatively curved setting, a landmark example was provided by Gromoll and Meyer (see \cite{gromoll}), who showed that the exotic sphere $\Sigma^7_{2,-1}$ can be expressed as a biquotient $\Sp(2) \biq \Sp(1)$. This led J. Eschenburg (see \cite{eschenburg}) to initiate the systematic study of biquotients, which provided a new way to construct spaces of positive curvature. Further examples include an infinite family in dimension $13$ known as Bazaikin spaces (see \cite{Bazaikin}), as well as a $7$-dimensional space homeomorphic to the unit tangent bundle of the $4$-sphere $T_1 S^4$ (see \cite{dearricott, groveverdianiziller}). Apart from this single exception, all known examples of simply connected manifolds with positive curvature arise as quotients of Lie groups or biquotients. For a survey of this discussion, we refer the reader to \cite{zillerpositive}.

When the isometric actions underlying such quotient constructions fail to be free, the resulting orbit spaces are no longer Riemannian manifolds. Nevertheless, under appropriate conditions, they naturally inherit the structure of Riemannian orbifolds. An important way to construct Riemannian orbifolds is from infinitesimally polar actions, i.e., isometric actions whose slice representations are polar. Recall that a proper and isometric action of a Lie group on a complete Riemannian manifold is polar if there exists a connected complete isometrically immersed submanifold meeting all orbits orthogonally. Indeed, the Lytchak-Thorbergsson theorem (see \cite{lytthor}), which establishes that an orbit space is a Riemannian orbifold if and only if the corresponding group action is infinitesimally polar, reduces the classification of orbifold biquotients to that of infinitesimally polar biquotient actions. In this direction, C. Gorodski and A. Lytchak classified infinitesimally polar actions of connected compact Lie groups on Euclidean spheres (see \cite{gorodskilyt}), and C. Gorodski and A. Kollross subsequently completed the classification on all compact rank-one symmetric spaces (see \cite{gorodskikollross}).

In this paper, we achieve a complete classification of the infinitesimally polar actions on the Lie group $\Sp(2)$ equipped with a bi-invariant Riemannian metric. It is worth noting that a special case of this problem was previously addressed by A. Gomes (see \cite{Andregomesthesis}), who classified the orbifold biquotients of the specific form $\Sp(2) \biq \Sp(1)$. Our approach consists of an exhaustive case-by-case analysis, determining all closed connected subgroups of $\Sp(2)$ up to conjugacy and examining all possible product, diagonal, and asymmetric subgroups of $\Sp(2) \times \Sp(2)$ that yield an orbifold biquotient of $\Sp(2)$ -- where a closed connected subgroup $U \subset \Sp(2) \times \Sp(2)$ is called \textit{asymmetric} if it is neither a product $H \times K$ nor conjugate in $\Sp(2) \times \Sp(2)$ to a diagonal subgroup $\Delta H = \{(h,h) : h \in H\}$ (see Section~\ref{Sec:subgroupsSp2}). More precisely, our main theorem is the following:

\begin{mainthm}\label{thm:main_classification}
Let $U \subset \Sp(2) \times \Sp(2)$ be a closed connected subgroup acting non-transitively on $\Sp(2)$ equipped with a bi-invariant Riemannian metric. The biquotient $\Sp(2) \biq U$ is a Riemannian orbifold -- equivalently, the $U$-action on $\Sp(2)$ is infinitesimally polar -- if and only if, up to conjugacy and permutation of the factors of $\Sp(2) \times \Sp(2)$, the action falls into one of the following classes:
\begin{enumerate}[label=\textnormal{(\roman*)}, itemsep=0.15cm]
    \item \label{thm:main_classificationitemi} Hermann actions: the action of conjugation by $U = \Delta \Sp(2)$ (of cohomogeneity $2$), or $U = H \times K$ with $H, K \in \{\Un(2), \SU(2) \times \SU(2)\}$ (of cohomogeneity $1$ or $2$);

    \item \label{thm:main_classificationitemii} Polar actions on $\Sp(2)/K$: the subgroup $U$ induces a polar action on $\CP^3$, $S^4$, or $S^7$ as listed in Table~\ref{table:polar_Sp2_K};
    
    \item \label{thm:main_classificationitemiii} Non-polar and non-almost-free actions: $U = H \times K$ induces a non-polar, infinitesimally polar action $H \curvearrowright \Sp(2)/K$, which is one of the following three cases:
    \begin{enumerate}[label=\textnormal{(\alph*)}]
        \item $\Un(2) \times (\SU(2) \times \Un(1))$ acting as $\Un(2) \curvearrowright \CP^3$, of cohomogeneity $2$;
        \item $\SU(2)_\Delta \times (\SU(2) \times \Un(1))$ acting as $\SU(2)_\Delta \curvearrowright \CP^3$, of cohomogeneity $3$;
        \item $\Un(2) \times \SU(2)_{\low}$ acting as $\Un(2) \curvearrowright S^7$, of cohomogeneity $3$;
    \end{enumerate}

    \item \label{thm:main_classificationitemiv} Almost-free actions: $U$ acts almost freely on $\Sp(2)$, and is either $\{e\} \times K$ for a proper closed connected subgroup $K \subsetneq \Sp(2)$ (yielding the homogeneous space $\Sp(2)/K$), or a product subgroup $H \times K$ or an asymmetric subgroup, as listed in Table~\ref{table:almost_free}.
\end{enumerate}
\end{mainthm}

\vspace{0.2cm}
The full classification presented in Theorem~\ref{thm:main_classification} combines classical results with new contributions. Hyperpolar actions on compact Lie groups -- comprising the Hermann actions and cohomogeneity one actions -- were classified by Kollross~\cite{kollross}. By exploiting the natural isometry between the orbit spaces of $H \times K \curvearrowright \Sp(2)$ and $H \curvearrowright \Sp(2)/K$, together with the classifications of infinitesimally polar actions on spheres due to Gorodski and Lytchak~\cite{gorodskilyt} and on complex projective spaces by Gorodski and Kollross~\cite{gorodskikollross} as well as the classification of polar actions on complex projective spaces due to Podestà and Thorbergsson~\cite{podesta-thorbergsson}, we determine all actions that are globally polar on $\Sp(2)/K$ (case~\ref{thm:main_classificationitemii}) as well as the three exceptional cases that are infinitesimally polar without being polar or almost free (case~\ref{thm:main_classificationitemiii}). Finally, our main contribution completes the classification for $\Sp(2)$ by explicitly determining all almost-free actions -- both of product and asymmetric type (case~\ref{thm:main_classificationitemiv}) -- as well as the asymmetric subgroups inducing polar actions on $S^7$ (in case~\ref{thm:main_classificationitemii}), and ruling out all remaining asymmetric subgroups. Moreover, since $\SO(5) = \Sp(2)/\{\pm I_4\}$ has trivial center, the identity component of its isometry group is $\Isom_0(\SO(5)) = \SO(5) \times \SO(5)$, and Theorem~\ref{thm:main_classification} immediately extends to the complete classification of infinitesimally polar actions on $\SO(5)$, as established in Subsection~\ref{subsec:classification_so(5)}:

\begin{maincor}\label{cor:SO5_classification}
Let $\psi \colon \Sp(2) \to \SO(5)$ be the universal double covering and set $\Psi \coloneqq \psi \times \psi \colon \Sp(2) \times \Sp(2) \to \SO(5) \times \SO(5)$. A closed connected subgroup $\tilde{U} \subset \SO(5) \times \SO(5)$ acting non-transitively on $\SO(5)$ equipped with a bi-invariant Riemannian metric acts infinitesimally polarly -- equivalently, the biquotient $\SO(5) \biq \tilde{U}$ is a Riemannian orbifold -- if and only if, up to conjugacy and permutation of the factors of $\SO(5) \times \SO(5)$, $\tilde{U} = \Psi(U)$ for one of the subgroups $U \subset \Sp(2) \times \Sp(2)$ classified in Theorem~\ref{thm:main_classification}.
\end{maincor}

In Section~\ref{Sec:inf_polar}, we establish Theorem~\ref{thm:main_classification} through a systematic case-by-case analysis. First, in Subsection~\ref{subsec:U2_subgroups}, we classify the infinitesimally polar actions of product subgroups $H \times K \subset \Un(2) \times \Un(2)$ with $H \subsetneq \Un(2)$. By analyzing the slice representations, we show that such an action is infinitesimally polar if and only if it is almost free. Furthermore, product actions $H \times K$ on $\Sp(2)$ that induce actions on the compact rank-one symmetric spaces $\CP^3$, $S^4$, $S^7$, and on the Berger space $B^7$ are analyzed in Subsections~\ref{subsec:sym-spaces}--\ref{section:berger-space}, whereas diagonal and asymmetric subgroup actions are classified in Subsections~\ref{seccao:diagonal} and~\ref{subsec:asymmetric}. The definitions of the subgroups in Table~\ref{table:polar_Sp2_K} and in Table~\ref{table:almost_free} can be found in either Section~\ref{Sec:subgroupsSp2} or Subsection~\ref{subsec:asymmetric}.

As a geometric application of Theorem~\ref{thm:main_classification}, we turn to the positively curved Berger space $B^7 = \Sp(2)/\SU(2)_{\irr}$ (see \cite{berger}) and its orbifold quotients obtained via O'Neill's submersion formula \cite{oneill}. While the almost-free action of $\SU(2)_{\low}$ on $B^7$ recovers the known weighted quaternionic projective space $\HP^1[2]$, our classification yields new examples of positively curved Riemannian orbifolds with non-trivial singular locus, which we prove in Section~\ref{sec:pos_curv}. First, we obtain an infinite family in dimension $6$, for which we define a parametrization for the one-dimensional tori of $\Sp(2)$ as $T^1_{m,n} \coloneqq \{ \diag(z^m, z^n, z^{-m},z^{-n}): z \in \Un(1) \}$, where $m,n \in \Z$ are coprime.  

\begin{maincor}\label{cor:positive_curvature}
For every pair of coprime integers $(m,n)$ such that $\{|m|, |n|\} \neq \{1,3\}$, the biquotient
\[
\mathcal{O}^6_{m,n} \coloneqq \Sp(2) \biq \big(T^1_{m,n} \times \SU(2)_{\irr}\big)
\]
is a $6$-dimensional positively curved Riemannian orbifold with non-trivial singular locus. 
\end{maincor}

In dimension $4$, the almost-free action of $\SU(2)_{\Delta}$ on $B^7$ produces a positively curved orbifold that is fundamentally distinct from $\HP^1[2]$ (whose singular locus has cyclic isotropy $\Z_3$), exhibiting non-cyclic isotropy and failing to be globally covered by a smooth manifold (i.e., a \textit{bad} orbifold), while its underlying topological space remains homeomorphic to $S^4$:

\begin{maincor}\label{cor:positive_curvature_4d}
The biquotient
\[\mathcal{Q}^4 \coloneqq \Sp(2) \biq \big(\SU(2)_{\Delta} \times \SU(2)_{\irr}\big)\]
is a $4$-dimensional positively curved bad Riemannian orbifold whose underlying topological space is homeomorphic to $S^4$. Its singular locus is homeomorphic to the space $D^2 \cup_{z \mapsto z^3} S^1$ obtained by attaching $D^2$ to $S^1$ via $z \mapsto z^3$, with local isotropy group $\Z_2$ on the open disk and $\Z_2 \times \Z_2$ on the circle $S^1$.
\end{maincor}

The proof of Corollary~\ref{cor:positive_curvature_4d} relies on realizing $\mathcal{Q}^4$ via $\SO(5)$ as an $\SO(3)$-quotient of the Stiefel manifold $V_2(\Sym_0(3,\R))$, together with an isometric circle action that preserves the singular locus and rotates the disk $D^2$.

\begin{table}[htb!]
    \centering
    \renewcommand{\arraystretch}{1.5}
    \begin{tabular}{|c|c|c|}
        \hline
        $U$ & \textsc{Action} & \textsc{Cohomogeneity} \\ \hline \hline
        $(\SU(2) \times \SU(2)) \times (\SU(2) \times \Un(1))$ & $\SU(2) \times \SU(2) \curvearrowright \CP^3$ & $1$ \\ \hline 
        $(\SU(2) \times \Un(1)) \times (\SU(2) \times \Un(1))$ & $\SU(2) \times \Un(1) \curvearrowright \CP^3$ & $2$ \\ \hline
        $\SU(2)_{\low} \times (\SU(2) \times \SU(2))$ & $\SU(2)_{\low} \curvearrowright S^4$ & $1$ \\ \hline
        $\SU(2)_{\Delta} \times (\SU(2) \times \SU(2))$ & $\SU(2)_{\Delta} \curvearrowright S^4$ & $2$ \\ \hline
        $\SU(2)_{\irr} \times (\SU(2) \times \SU(2))$ & $\SU(2)_{\irr} \curvearrowright S^4$ & $1$ \\ \hline
        $T^2 \times (\SU(2) \times \SU(2))$ & $T^2 \curvearrowright S^4$ & $2$ \\ \hline 
        $T^1_{m,n} \times (\SU(2) \times \SU(2)), \ m=\pm n$ & $T^{1}_{m,n} \curvearrowright S^4$ & $3$ \\ \hline
        $\SU(2)_{\low} \times \SU(2)_{\low}$ & $\SU(2)_{\low} \curvearrowright S^7$ & $4$ \\ \hline
        $\Delta\SU(2)_{\mathrm{up}} \cdot (\SU(2)_{\low} \times \SU(2)_{\low})$ & $\SU(2) \times \SU(2) \curvearrowright S^7$ & $2$ \\ \hline
        $U_{m,k}, \ m=\pm k$ & $\SU(2) \times \Un(1) \curvearrowright S^7$ & $3$ \\ \hline
    \end{tabular}
    \vspace{0.3cm}
    \caption{Subgroups $U \subset \Sp(2) \times \Sp(2)$ inducing polar actions on $\Sp(2)/K$.}
    \label{table:polar_Sp2_K}
\end{table}

\begin{table}[htb!]
    \centering
    \renewcommand{\arraystretch}{1.5}
    \begin{tabular}{|c|c|c|}
        \hline 
        $U$ & \textsc{Condition} & \textsc{$\dim (\Sp(2) \biq U)$} \\ \hline \hline
        $\SU(2)_\Delta \times T^1_{m,n}$ & $|m| \neq |n|$ & $6$ \\ \hline
        $T^1_{m,n} \times T^1_{a,b}$ & $\{ \vert a \vert, \vert b \vert \} \neq \{ \vert m \vert, \vert n \vert \}$ & $8$ \\ \hline
        $\SU(2)_{\Delta} \times \SU(2)_{\low}$ & --- & $4$ \\ \hline
        $\SU(2)_{\irr} \times \SU(2)_{\low}$ & --- & $4$ \\ \hline
        $T^1_{m,n} \times \SU(2)_{\low}$ & $m \cdot n \neq 0$ & $6$ \\ \hline
        $T^1_{m,n} \times \SU(2)_{\irr}$ & $\{|m|,|n|\} \neq \{1,3\}$ & $6$ \\ \hline
        $\SU(2)_\Delta \times \SU(2)_{\irr}$ & --- & $4$ \\ \hline
        $\Delta^{\theta_c}(\Un(2), \SU(2) \times \Un(1))$ & $|c| \neq \frac{3}{2}$ & $6$ \\ \hline
        $U^c_{m,n}$ & $c(\pm m \pm n) \neq 2$ & $6$ \\ \hline
        $\Un(1) \cdot \Delta^{\theta}(\SU(2)_{\Delta}, \SU(2)_{\low})$ & --- & $6$ \\ \hline
        $\Delta^{\theta}(\SU(2)_{\irr}, \SU(2)_{\Delta})$ & --- & $7$ \\ \hline
        $\Delta^{\theta}(\SU(2)_{\irr}, \SU(2)_{\low})$ & --- & $7$ \\ \hline
        $\Delta^{\theta}(\SU(2)_{\low}, \SU(2)_{\Delta})$ & --- & $7$ \\ \hline
        $V^c_{m,n}$ & $|c| \notin \{|m|, |n|\}$ & $6$ \\ \hline
        $\Delta^{\theta} T^2$ & $\det(\id - sA) \neq 0, \ \forall s \in D_4$ & $8$ \\ \hline
        $T^2_a(k_1,k_2,m,n)$ & $a \langle (-k_2, k_1), s(m,n) \rangle \neq 1, \ \forall s \in D_4$ & $8$ \\ \hline
        $U^{r,s}_{m,n}$ & $\{|m|, |n|\} \neq \{|r|, |s|\}$ & $9$ \\ \hline
        $H_{p,q}(m,n,r,s)$ & $p, q \neq 0, \ |p| \neq |q|, \ \gcd(p,q) = 1$ & $9$ \\ \hline
    \end{tabular}
    \vspace{0.3cm}
    \caption{Closed connected subgroups $U \subset \Sp(2) \times \Sp(2)$ (other than $\{e\} \times K$) acting almost freely on $\Sp(2)$.}
    \label{table:almost_free}
\end{table}

\medskip
\textit{Acknowledgments.} We would like to thank Claudio Gorodski for his support and several valuable discussions related to this work. N. Caballero has been supported by Grant 2025/00282-3 from the S\~ao Paulo Research Foundation (FAPESP). D. Rotmeister has been supported by Grant 2024/01929-8 from the S\~ao Paulo Research Foundation (FAPESP).

\section{Preliminaries}

In this section, we review the fundamental properties of biquotients and (infinitesimally) polar actions. For a comprehensive treatment of these topics, we refer the reader to \cite{gorodski, berndt-console-olmos, alex-bettiol}. 

\subsection{Biquotients}

Let $G$ be a compact Lie group endowed with a bi-invariant metric, and consider a closed subgroup $U$ of $G \times G$ acting on $G$ by $(u_1,u_2) \cdot g \coloneqq u_1 g u_2^{-1}$ for all $(u_1,u_2) \in U$ and $g \in G$. The corresponding orbit space is called a \textit{biquotient} and is denoted by $G \biq U$. Endowing its principal stratum with the quotient metric so that the canonical projection $G \to G \biq U$ restricts to a Riemannian submersion on the union of principal orbits, O'Neill's formula \cite{oneill} yields non-negative sectional curvature on this open dense stratum. Whenever $G \biq U$ is a Riemannian orbifold, this non-negative curvature extends to the entire orbifold by continuity on local charts. 

\begin{rmk}
    It is worth noting that biquotients are generalizations of homogeneous spaces. Indeed, a homogeneous space $G/H$ can be written as $G \biq U$ where $U = \{e\} \times H$. Furthermore, if a closed subgroup $K \subset G$ acts on $G/H$, then the quotient $K \backslash G /H$ is a biquotient with $U = K \times H$. 
\end{rmk}

\begin{rmk}
    Let $U \subset G \times G$ be a subgroup. For every $a,b \in G$, consider 
    \[
    U_{a,b} = \{ (au_1a^{-1}, bu_2 b^{-1}) : (u_1,u_2) \in U \} .
    \]
    The $U$- and $U_{a,b}$-actions are conjugate via the isometry $g \in G \mapsto agb^{-1} \in G$. Moreover, we will not distinguish the order of the two factors of $G \times G$ since the actions of $U$ and $\{ (u_2, u_1) : (u_1, u_2) \in U \}$ on $G$ are conjugate by the isometry $g \mapsto g^{-1}$.
\end{rmk}

Following \cite{lytthor}, we now recall the definition of a Riemannian orbifold, which serves as the main object of study in this paper.

\begin{Def}\label{Def:Riemorb}
    An intrinsic metric space $X$ is called a \textit{Riemannian orbifold} if every $x \in X$ admits a neighborhood $V$ isometric to a quotient $M / \Gamma$, where $M$ is a connected Riemannian manifold and $\Gamma \subset \Isom(M)$ is a finite subgroup of isometries.
\end{Def}

Let $X$ be a Riemannian orbifold. For a given point $x \in X$, consider a local neighborhood isometric to a quotient $M/\Gamma$ and let $p \in M$ be a lift of $x$ under the natural projection $\pi^{\mathrm{orb}} \colon M \to M/\Gamma$. The isotropy subgroup at $x$ is defined by $\Gamma_p \coloneqq \{ g \in \Gamma : g \cdot p = p \}$. Up to isomorphism, $\Gamma_p$ depends only on the point $x$ and neither on the local chart $M/\Gamma$ nor on the choice of the lift $p$. The point $x$ is said to be \textit{regular} if $\Gamma_p$ is trivial and \textit{singular} otherwise. The set of all singular points is called the \textit{singular locus} of $X$, denoted by $\Sigma(X)$. Moreover, a Riemannian orbifold $X$ is called \textit{good} if it is globally isometric to a quotient $M/\Gamma$ of a Riemannian manifold $M$ by a discrete group $\Gamma$ of isometries acting properly, and \textit{bad} otherwise (see \cite{caramello}).

An action $U \curvearrowright G$ is called \textit{almost free} (resp. \textit{free}) if all the isotropy subgroups are finite (resp. trivial). Under these conditions, the biquotient $G \biq U$ is a Riemannian orbifold (resp. a smooth manifold), see \cite[Proposition~1.5.1]{caramello}.

\begin{ex} 
The \textit{weighted quaternionic projective space} $\mathbb{H} P^n[d_1, \ldots, d_k]$, with $d_i > 0$ and $\sum d_i = n+1$, is defined as the quotient of the sphere $S^{4n+3} \subset \mathbb{H}^{n+1} = \mathbb{H}^{d_1} \oplus \cdots \oplus \mathbb{H}^{d_k}$ by the left action of $S^3 = \Sp(1)$ given by $\rho = \rho_{d_1} \oplus \dots \oplus \rho_{d_k}$, where $\rho_d$ is the irreducible $d$-dimensional quaternionic representation ($4d$-dimensional real representation) of $\Sp(1)$. The simplest example of a non-trivial weighted quaternionic projective space is $\mathbb{H}P^1[2]$. This orbifold can be obtained as the orbit space $\Sp(1)_{\low} \backslash \Sp(2) /\Sp(1)_{\irr}$, where the group on the right corresponds to the irreducible $2$-dimensional quaternionic representation, and the group on the left corresponds to the embedding $\{ \diag(1,q)  :  q \in \Sp(1) \} \subset \Sp(2)$. 
\end{ex}

\subsection{Generalities on Polar Actions}

Let $M$ be a complete Riemannian manifold and consider an isometric action $G \curvearrowright M$ of a compact Lie group $G$. The \textit{isotropy representation} of the action at $x \in M$ is the linear action of the isotropy subgroup $G_x$ on the tangent space $T_x M$ obtained by differentiating the $G_x$-action at $x$. It is straightforward to verify that $T_x (G \cdot x)$ and the normal space $\nu_x(G \cdot x) =(T_x (G \cdot x))^\perp$ to the orbit through $x$ are invariant subspaces of the isotropy representation. The restriction of the isotropy representation to $\nu_x(G \cdot x)$ is called the \textit{slice representation} at $x$. Recall that the action $G \curvearrowright M$ is called \textit{polar} if there exists a complete connected immersed submanifold $\Sigma$ of $M$, called a \textit{section}, which intersects all the orbits and such that $\Sigma$ is perpendicular to every orbit it meets. In addition, whenever a section $\Sigma$ is flat with respect to the induced Riemannian metric, the action is called \textit{hyperpolar}.

\begin{Def}
 An action $G \curvearrowright M$ is called \textit{infinitesimally polar} if the slice representation at $p \in M$ is polar for all $p \in M$. 
\end{Def}

Infinitesimal polarity serves as the precise local criterion governing the global differential structure and singularities of the quotient space $M/G$. This connection was established by A. Lytchak and G. Thorbergsson.

\begin{thm}[\bfseries Lytchak-Thorbergsson \cite{lytthor}]\label{Thm:LTequiv}
    Let $G \curvearrowright M$ be an isometric action of a compact Lie group on a Riemannian manifold. The orbit space $M/G$ is a Riemannian orbifold if and only if $G \curvearrowright M$ is infinitesimally polar.
\end{thm}

Before proceeding, we state a celebrated structural result concerning linear polar actions. In 1985, J. Dadok completely classified polar representations by proving that, up to orbit equivalence, they arise as the linear isotropy representations of symmetric spaces. More precisely:

\begin{thm}[\bfseries Dadok \cite{dadok}] \label{Thm:Dadok} Let $G$ be a closed, connected subgroup of the special orthogonal group $\SO(V)$ of a Euclidean space $V$, and assume that its action on $V$ is polar. Then the following hold: 
\begin{enumerate}[label=\textnormal{(\roman*)}]
    \item \label{Thm:Dadok(i)} If $V$ decomposes as the sum of $G$-invariant subspaces $V_1 \oplus V_2$ then the $G$-action on each $V_i$ is polar for $i \in \{1,2\}$ and the cohomogeneity of the $G$-action on $V$ is the sum of the cohomogeneities of the $G$-actions on $V_i$, $i \in \{1,2\}$;

    \item \label{Thm:Dadok(ii)} If $\hat{G}$ denotes the maximal connected subgroup of $\SO(V)$ having the same orbits as $G$, then the action of $\hat{G}$ on $V$ is conjugate to the isotropy representation of a symmetric space $M$. In particular, the rank of $M$ equals the cohomogeneity of the $G$-action on $V$, i.e., $\rank(M) = \cohom(G, V)$.
\end{enumerate}
\end{thm}

\begin{cor}\label{corolario:toro}
Let $T^r$ be an $r$-dimensional torus acting almost effectively and unitarily on $\C^s$ without non-trivial fixed vectors. If $r \ne s$, then the action is not polar. In particular, if $r < s$, any unitary action of $T^r$ on $\C^s$ without non-trivial fixed vectors is not polar.
\end{cor}

\begin{proof}
Since every irreducible complex representation of a torus is one-dimensional and there are no non-trivial fixed vectors, the complex space $\C^s$ decomposes into orthogonal weight spaces as $\C^s = \C_{\omega_1} \oplus \dots \oplus \C_{\omega_s}$ with $\omega_i \neq 0$ for all $i$. The cohomogeneity of the action on each factor is $1$, and therefore the sum of the cohomogeneities is $s$. On the other hand, since $T^r$ is Abelian and acts almost effectively on $\C^s$, the stabilizer of a generic vector $(v_1, \dots, v_s)$ (with $v_i \neq 0$ for all $i$) coincides with the finite ineffective kernel, so the principal orbits have dimension $r$ and the cohomogeneity on $\C^s$ is $2s - r$. Hence, if $r \neq s$, then $2s - r \neq s$, and by Theorem~\ref{Thm:Dadok}\ref{Thm:Dadok(i)} the action cannot be polar. Finally, if $r < s$ without assuming almost effectiveness, passing to the effective quotient $T^{r'}$ yields $r' \leq r < s$, so the same conclusion follows.
\end{proof}

In addition to the previous results, it is important to note that a linear representation of a compact Lie group is polar if and only if the restriction to its identity component is polar. In order to identify almost-free actions on $\Sp(2)$, we state the following lemma.

\begin{lemma}\label{lemma:almostfree}
Let $G$ be a compact Lie group, $U \subset G \times G$ a closed connected subgroup, and $\mathfrak{t}_u \subset \mathfrak{u} = \mathrm{Lie}(U)$ a maximal Abelian subalgebra. Then, the $U$-action on $G$ is almost free if and only if for all non-zero $(X_1,X_2) \in \mathfrak{t}_u \subset \mathfrak{u} \subset \mathfrak{g} \oplus \mathfrak{g}$ we have that $X_1 - \mathrm{Ad}_g(X_2) \neq 0$ for all $g \in G$. 
\end{lemma}
\begin{proof}
By definition, the $U$-action on $G$ is almost free if and only if all isotropy groups are finite, which is equivalent to $\mathrm{Lie}(U_g)=0$ for every $g \in G$, where $U_g$ denotes the isotropy group at $g \in G$. Recalling that the isotropy subgroup at a point $g \in G$ is given by $U_g = \{ (h,k) \in U: hgk^{-1}=g \}$, it turns out that its Lie algebra is given by 
$$\mathrm{Lie}(U_g) = \{(X_1,X_2) \in \mathfrak{u}: X_1 = \mathrm{Ad}_g(X_2) \}.$$
Hence, the action is almost free if and only if $X_1 - \mathrm{Ad}_g(X_2) \neq 0$ for all non-zero $(X_1, X_2) \in \mathfrak{u}$. Finally, since every element of the Lie algebra $\mathfrak{u}$ is conjugate to an element of a maximal Abelian subalgebra $\mathfrak{t}_u$, it suffices to check the condition for $(X_1,X_2)\in\mathfrak{t}_u$.
\end{proof}

\section{Subgroups of $\Sp(2) \times \Sp(2)$}\label{Sec:subgroupsSp2}
In this section, we present our method to identify all the closed connected subgroups $U$ of the product $\Sp(2) \times \Sp(2)$ and hence to classify the biquotients $\Sp(2) \biq U$. To this end, we review the standard setup needed to introduce Goursat's Lemma.  

Let $\mathfrak{g}_1$ and $\mathfrak{g}_2$ be Lie algebras. For each $i \in \{1,2\}$, denote by $\pr_i : \mathfrak{g}_1 \oplus \mathfrak{g}_2 \to \mathfrak{g}_i$ the canonical projection. Given a subalgebra $\mathfrak{u} \subset \mathfrak{g}_1 \oplus \mathfrak{g}_2$, one considers $\mathfrak{h}_i \coloneqq \pr_i(\mathfrak{u})$, which is a subalgebra of $\mathfrak{g}_i$. Under the natural identifications $\mathfrak{g}_1 \oplus \{0\} \cong \mathfrak{g}_1$ and $\{0\} \oplus \mathfrak{g}_2 \cong \mathfrak{g}_2$, it is straightforward to check that $\mathfrak{n}_1 = \mathfrak{u} \cap (\mathfrak{g}_1 \oplus \{0\})$ and $\mathfrak{n}_2 = \mathfrak{u} \cap (\{0\} \oplus \mathfrak{g}_2)$ are ideals of $\mathfrak{h}_1$ and $\mathfrak{h}_2$, respectively. Hence, the quotient space $\mathfrak{h}_i/\mathfrak{n}_i$ is a Lie algebra and the map $\theta : \mathfrak{h}_1/\mathfrak{n}_1 \to \mathfrak{h}_2/\mathfrak{n}_2$ given by $\theta(X + \mathfrak{n}_1) = Y + \mathfrak{n}_2$ for any $(X,Y) \in \mathfrak{u}$ is an isomorphism. Conversely, given a quintuple $\{\mathfrak{h}_1, \mathfrak{n}_1 , \mathfrak{h}_2, \mathfrak{n}_2, \theta \}$, the subalgebra $\mathfrak{u} \subset \mathfrak{g}_1 \oplus \mathfrak{g}_2$ is defined by the fiber product $\mathfrak{h}_1 \times _{\mathfrak{h}_2/\mathfrak{n}_2} \mathfrak{h}_2$, i.e., by the commutative diagram
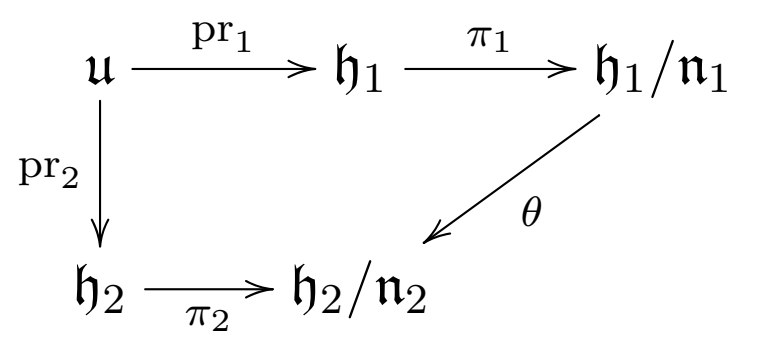
\begin{align*}
     \xymatrix{
     \mathfrak{u} \ar[r]^-{\pr_1} \ar[d]_-{\pr_2} & \mathfrak{h}_1 \ar[r]^-{\pi_1} & \mathfrak{h}_1/\mathfrak{n}_1 \ar[ld]^{\theta} \\
     \mathfrak{h}_2 \ar[r]_-{\pi_2} & \mathfrak{h}_2/\mathfrak{n}_2 &
     }
\end{align*}
This implies that 
\begin{align}\label{Goursatslemma:subalgu}
    \mathfrak{u} & = \{ (X,Y) \in \mathfrak{h}_1 \oplus \mathfrak{h}_2 : \theta(\pi_1(X)) = \pi_2(Y) \} \\
    &= \{ (X_1+X_2, Y_1+\theta(X_2)) : X_1 \in \mathfrak{n}_1, \ X_2 \in \mathfrak{h}_1/\mathfrak{n}_1, \ Y_1 \in \mathfrak{n}_2 \}. \nonumber
\end{align}
These two assignments are mutually inverse, establishing the following result, which can also be found in \cite{chapov-kov}:

\begin{thm}[\bfseries Goursat's Lemma for Lie Algebras \cite{goursat}] 
There is a bijective correspondence between Lie subalgebras $\mathfrak{u}$ of $\mathfrak{g}_1 \oplus \mathfrak{g}_2$ and quintuples $\{\mathfrak{h}_1, \mathfrak{n}_1, \mathfrak{h}_2, \mathfrak{n}_2, \theta\}$, where $\mathfrak{h}_i \subseteq \mathfrak{g}_i$ is a Lie subalgebra, $\mathfrak{n}_i \unlhd \mathfrak{h}_i$ is an ideal for $i\in\{1,2\}$, and $\theta: \mathfrak{h}_1/\mathfrak{n}_1 \to \mathfrak{h}_2/\mathfrak{n}_2$ is a Lie algebra isomorphism. 
\end{thm}

To exhaust all possibilities for closed connected subgroups $U \subset \Sp(2) \times \Sp(2)$, we apply Goursat's Lemma as follows. Fix a Lie subalgebra $\mathfrak{h}_1 \subset \mathfrak{sp}(2)$ and determine its ideals $\mathfrak{n}_1$, as well as the corresponding quotients $\mathfrak{h}_1/\mathfrak{n}_1$. Then, scan through all Lie subalgebras $\mathfrak{h}_2$ and their ideals $\mathfrak{n}_2$ such that $\mathfrak{h}_2/\mathfrak{n}_2$ is isomorphic to $\mathfrak{h}_1/\mathfrak{n}_1$. Finally, up to conjugation, determine all isomorphisms $\theta: \mathfrak{h}_1/\mathfrak{n}_1 \to \mathfrak{h}_2/\mathfrak{n}_2$ such that the corresponding quintuple gives rise to a Lie subalgebra $\mathfrak{u}$ which integrates to a closed connected subgroup $U \subset \Sp(2) \times \Sp(2)$.

Applying the method above, a closed connected subgroup $U$ of $\Sp(2) \times \Sp(2)$ has one of the following forms, up to conjugation and permutation of the factors: $U$ is a product $H \times K$ of closed connected subgroups of $\Sp(2)$; $U$ is a \textit{twisted diagonal} $\Delta^\sigma H = \{ (h,\sigma(h)): h \in H \}$, where $\sigma \in \mathrm{Aut}(\Sp(2))$; or $U$ is an \textit{asymmetric subgroup}, that is, $U$ is such that the Goursat's Lemma quintuple $\{\mathfrak{h}_1, \mathfrak{n}_1, \mathfrak{h}_2, \mathfrak{n}_2, \theta \}$ has a non-trivial isomorphism $\theta$ and $U$ is neither a twisted diagonal nor a product subgroup. In fact, if the quintuple is such that $\mathfrak{n}_1 = \mathfrak{h}_1$ and $\mathfrak{n}_2 = \mathfrak{h}_2$, then the subgroup $U$ is a product $H_1 \times H_2$, where $\mathfrak{h}_i = \mathrm{Lie}(H_i)$, $i = 1,2$. If $\mathfrak{h}_1 = \mathfrak{h}_2$ (up to conjugation in $\Sp(2)$), $\mathfrak{n}_1 = \mathfrak{n}_2 = \{ 0 \}$, and $\theta$ is the restriction of $\d\sigma_e$ for some $\sigma \in \Aut(\Sp(2))$, then $U$ is conjugate to a twisted diagonal $\Delta^\sigma H$. Finally, if the quintuple is not in one of the cases above, then it yields an asymmetric subgroup $U$. It is now necessary to identify the Lie subalgebras of $\mathfrak{sp}(2)$. 

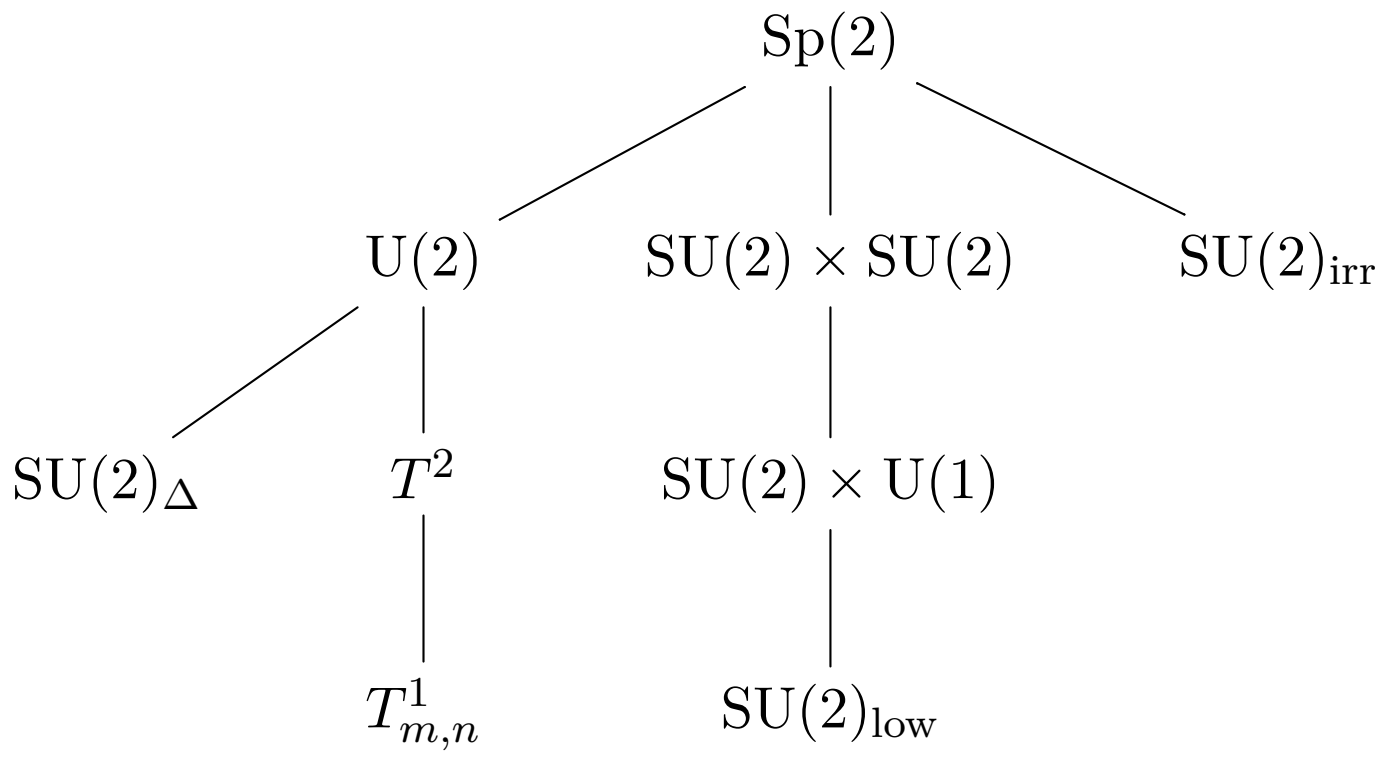
\begin{figure}[htb!]
\[
\xymatrix{
& & \Sp(2) \ar@{-}[ld] \ar@{-}[d] \ar@{-}[rd] & & \\
& \Un(2) \ar@{-}[ld] \ar@{-}[d] 
& \SU(2) \times \SU(2) \ar@{-}[d] 
& \SU(2)_{\mathrm{irr}} & \\
\SU(2)_{\mathrm{\Delta}} 
& T^2 \ar@{-}[d] 
& \SU(2) \times \Un(1)  \ar@{-}[d] & \\
& T_{m,n}^1 & \SU(2)_{\mathrm{low}} & &
}
\]
\caption{Closed connected subgroups of $\Sp(2)$ up to conjugacy.}
\label{Fig:treeSp2}
\end{figure}

We now review the systematic classification of the closed connected subgroups of $\Sp(2)$ due to Dynkin (see \cite{dynkin}), illustrated in Figure~\ref{Fig:treeSp2}, in order to determine the Lie subalgebras of $\mathfrak{sp}(2)$ that can be identified as the Lie algebra of a closed connected subgroup $H$ of $\Sp(2)$. At the top of the tree is $\Sp(2)$ and its branches represent the maximal subgroups, up to conjugation in $\Sp(2)$. From each of these maximal subgroups, we extend further branches to represent their own subgroups of equal or lower rank (omitting cross-inclusions between different branches for clarity).

Henceforth, we identify the quaternions $\mathbb{H}$ with $\C^2$ by viewing $\mathbb{H}$ as a right complex vector space. Specifically, any $q \in \mathbb{H}$ can be uniquely written in the form $q = z_1 + jz_2$, with $z_1, z_2 \in \C$. Under this convention, the group of unit quaternions is
\begin{align*}
    \Sp(1) \coloneqq \{ q \in \mathbb{H} : q \thickbar{q}=1 \} \approx S^3
\end{align*}
and its Lie algebra $\mathfrak{sp}(1)$ is canonically isomorphic to $\mathfrak{su}(2)$. Generalizing to the quaternionic plane $\mathbb{H}^2$, the \textit{symplectic group} is the compact Lie group of $\mathbb{H}$-linear isometries defined by
\[ \Sp(2) \coloneqq \{ A \in \GL(2, \mathbb{H}) : AA^\ast = I \}. \]
The inner product $\langle X,Y \rangle = \mathrm{Re} \tr(XY^\ast)$ on $\mathfrak{sp}(2)$ induces a bi-invariant Riemannian metric on $\Sp(2)$ and is a multiple of the Cartan-Killing form of $\mathfrak{sp}(2)$. Moreover, restricting the embedding $M(2, \mathbb{H}) \hookrightarrow M(4, \C)$ given by
\begin{align*}
    A+jB \longmapsto \begin{pmatrix}
        A & -\thickbar{B} \\
         B & \thickbar{A}
    \end{pmatrix}
    , \qquad A,B \in M(2,\C),
\end{align*} 
we identify $\Sp(2)$ with a closed subgroup of $\SU(4)$. The maximal subgroups of $\Sp(2)$ are $\Un(2)$, $\SU(2) \times \SU(2)$, and $\SU(2)_{\mathrm{irr}}$. The first one is embedded as $\Un(2) = \left\{ \diag(A, \thickbar{A}) : A \in \Un(2) \right\} \subset \Sp(2)$, whereas $\SU(2) \times \SU(2)$ embeds as 
$$\SU(2) \times \SU(2) =  \left\{\begin{pmatrix} \alpha_1 & 0 & \beta_1 & 0 \\ 0 & \alpha_2 & 0 & \beta_2 \\ -\thickbar{\beta}_1 & 0 & \thickbar{\alpha}_1 & 0 \\ 0 & -\thickbar{\beta}_2 & 0 & \thickbar{\alpha}_2 \end{pmatrix}: \alpha_i, \beta_i \in \C, \vert \alpha_i \vert^2 + \vert \beta_i \vert^2 = 1 \right\}.$$
The maximal subgroup $\SU(2) \times \SU(2)$ of $\Sp(2)$ contains the $4$-dimensional subgroup $\SU(2) \times \Un(1)$ given by
\[\SU(2) \times \Un(1) = \left\{\begin{pmatrix} \alpha & 0 & \beta & 0 \\ 0 & z & 0 & 0 \\ -\thickbar{\beta} & 0 & \thickbar{\alpha} & 0 \\ 0 & 0 & 0 & \thickbar{z} \end{pmatrix}: \alpha, \beta, z \in \C, \vert \alpha \vert^2 + \vert \beta \vert^2 = 1, \vert z \vert = 1 \right\},\]
whereas the maximal subgroup $\SU(2)_{\mathrm{irr}}$ of $\Sp(2)$ corresponds to the irreducible $4$-dimensional complex representation of $\SU(2)$ which embeds as 
\begin{align}\label{Def:SU(2)irr}
    \SU(2)_{\mathrm{irr}} \coloneqq \left\{ \begin{pmatrix} \alpha^3 & \sqrt{3} \alpha \thickbar{\beta}^2 & -\thickbar{\beta}^3 & -\sqrt{3}\alpha^2 \thickbar{\beta} \\
\sqrt{3}\alpha \beta^2 & \alpha\thickbar{\alpha}^2 - 2\thickbar{\alpha}\beta \thickbar{\beta} & -\sqrt{3} \thickbar{\alpha}^2 \thickbar{\beta} & 2\alpha \thickbar{\alpha} \beta - \beta^2 \thickbar{\beta} \\ \beta^3 & \sqrt{3}\thickbar{\alpha}^2 \beta & \thickbar{\alpha}^3 & \sqrt{3} \thickbar{\alpha} \beta^2 \\
\sqrt{3}\alpha^2 \beta & \beta \thickbar{\beta}^2 - 2\alpha\thickbar{\alpha} \thickbar{\beta} & \sqrt{3} \thickbar{\alpha} \thickbar{\beta}^2 & \alpha^2 \thickbar{\alpha} - 2 \alpha \beta \thickbar{\beta}
\end{pmatrix}: \alpha, \beta \in \C, \vert \alpha \vert^2 + \vert \beta \vert^2 = 1  \right\} .
\end{align}

Now it is necessary to fix the notation for the toral subgroups, up to conjugation. For the one-dimensional tori, let 
\begin{align*}
    T^1_{m,n} & = \left\{ \diag(z^{m}, z^{n}, z^{-m}, z^{-n}) : z \in \Un(1) \right\} 
\end{align*}
where $m,n \in \Z$ with $\gcd(m, n) =1$. For the maximal torus of $\Sp(2)$, we consider its embedding as  
\begin{equation*}
    T^2 = \{ \diag(e^{i\theta_1}, e^{i\theta_2}, e^{-i\theta_1}, e^{-i\theta_2}) : \theta_1, \theta_2 \in \R \}.
\end{equation*}

In order to determine the non-maximal $\SU(2)$-subgroups of $\Sp(2)$, it suffices to analyze the $\SU(2)$-representations on $\C^4$ whose images are properly contained in $\Sp(2)$. The subgroup 
\begin{align}\label{Def:SU2low}
    \SU(2)_{\mathrm{low}} \coloneqq \left\{\begin{pmatrix} 1 & 0 & 0 & 0 \\ 0 & \alpha & 0 & \beta \\ 0 & 0 & 1 & 0 \\ 0 & -\thickbar{\beta} & 0 & \thickbar{\alpha} \end{pmatrix}: \alpha, \beta \in \C, \vert \alpha \vert^2 + \vert \beta \vert^2 = 1 \right\}
\end{align}
corresponds to the representation $\C \oplus \C \oplus \C^2$, where $\SU(2)$ acts by the standard representation on $\C^2$ and trivially on the two factors $\C$. On the other hand, the $\C^2 \oplus \thickbar{\C}^2$ representation of $\SU(2)$ characterizes the subgroup 
\begin{align}\label{Def:SUconj}
    \SU(2)_{\Delta} \coloneqq \{ \diag(A, \thickbar{A}): A \in \SU(2) \},
\end{align}
which is conjugate to the subgroup $\SU(2)$ that arises from its $\C^2 \oplus \C^2$ representation.

\section{The Classification}\label{Sec:inf_polar}

Our approach to identifying biquotients of $\Sp(2)$ that admit orbifold structures is based on the characterization provided in Theorem \ref{Thm:LTequiv}.  In what follows, we will explain our method. Let $H,K \subset \Sp(2)$ be closed connected subgroups. Consider the action $H \times K \curvearrowright \Sp(2)$ and its isotropy subgroup at the identity $e \in \Sp(2)$ given by
\begin{align}\label{HxK:isotropyate}
    (H \times K )_e = \{ (h,h) : h \in H \cap K \} = \Delta (H \cap K) \cong H \cap K.
\end{align}
The tangent and normal spaces to the orbit through $e$ are given, respectively, by
\begin{align*}
    T_e ((H \times K) \cdot e) = \mathfrak{h} + \mathfrak{k} \quad \textnormal{and} \quad \nu_e((H \times K) \cdot e) = \mathfrak{h}^\perp \cap \mathfrak{k}^\perp,
\end{align*}
where the orthogonal complements are taken with respect to the Cartan-Killing form of $\mathfrak{sp}(2) = \mathrm{Lie}(\Sp(2))$. It is easily checked that the isotropy representation at $e$ is equivalent to the adjoint representation of $\Sp(2)$ restricted to $H \cap K$. Therefore, the corresponding slice representation at $e$ is equivalent to the adjoint representation of $\Sp(2)$ restricted to $H \cap K$ on the invariant subspace $\mathfrak{h}^\perp \cap \mathfrak{k}^\perp$.

For a subgroup $U \subset \Sp(2) \times \Sp(2)$ acting on $\Sp(2)$ and a point $g \in \Sp(2)$, the isotropy subalgebra is $\mathfrak{u}_g = \{(X,Y) \in \mathfrak{u} : X = \Ad_g Y\}$. Pulling $T_g(U \cdot g)$ back to the identity via the isometry $\d R_{g^{-1}}$ yields 
$$V_g \coloneqq \{ X - \Ad_g Y : (X,Y) \in \mathfrak{u} \} \subset \mathfrak{sp}(2),$$ 
so the normal space $\nu_g(U \cdot g)$ is identified with $V_g^\perp \subset \mathfrak{sp}(2)$. In addition, the isotropy subgroup at some other point $g \in \Sp(2)$ with $g \neq e$ can be computed using \eqref{HxK:isotropyate}. Indeed, it is straightforward to check
\begin{align*}
    (H \times g Kg^{-1})_e = \{ (h, gkg^{-1}) \in H \times gKg^{-1} : hgk^{-1} =g \} \cong (H \times K)_g .
\end{align*}

Recall that all the automorphisms of $\Sp(2)$ are inner, i.e., for every $\sigma \in \Aut(\Sp(2))$, there exists $g \in \Sp(2)$ such that $\sigma(h) = ghg^{-1}$ for all $h \in \Sp(2)$. Consequently, for every closed connected subgroup $H \subset \Sp(2)$ and every $\sigma \in \Aut(\Sp(2))$, the twisted diagonal $\Delta^\sigma H$ is conjugate to $\Delta H$ in $\Sp(2) \times \Sp(2)$. Therefore, it remains to analyze the subgroups of the form $H \times K$ where $H$ and $K$ belong to the tree of Figure \ref{Fig:treeSp2}, the subgroups $\Delta H$ and the asymmetric subgroups given by Goursat's Lemma. It is worth pointing out that $\Un(2)$ and $\SU(2) \times \SU(2)$ are the only symmetric subgroups of $\Sp(2)$. Therefore, for any $H,K \in \{ \Un(2), \SU(2) \times \SU(2) \}$, the corresponding action of $H \times K$ on $\Sp(2)$ is a Hermann action (see \cite{kollross}).

For a closed connected subgroup $H \subset \Sp(2)$, we define the space
$$\mathrm{Fix}(H, \mathfrak{sp}(2)) \coloneqq \{ X \in \mathfrak{sp}(2): \Ad_g(X) = X, \ \text{for all} \ g \in H \}.$$
With this notation, we state the following result, which gives a condition to determine when an $H \times K$-action on $\Sp(2)$ can be infinitesimally polar given that the intersection $(H \cap K)^0$ is a torus. 

\begin{thm} Let $H,K$ be closed connected subgroups of $\Sp(2)$ such that $(H \cap K)^0 = T^r$, where $T^r$ is an $r$-dimensional torus and $r \in \{1,2\}$. If $\dim H + \dim K \leq 9 - \dim (\mathrm{Fix}(T^r,\mathfrak{sp}(2)))$, then the $H \times K$-action on $\Sp(2)$ is not infinitesimally polar.
\label{theo:intersecao-toro}
\end{thm}

\begin{proof} 
Let $V = \mathfrak{h}^\perp \cap \mathfrak{k}^\perp$ be the normal space to the orbit of $e \in \Sp(2)$ at this point. Denote by $V_0$ and $V_1$ the subspaces of $V$ where the isotropy subgroup $(H \times K)_e^0 \cong T^r$ acts trivially and non-trivially, respectively, by the restriction of $\mathrm{Ad} \coloneqq \mathrm{Ad}_{\Sp(2)}$ to $T^r$ on the space $V= V_0 \oplus V_1$. By Corollary~\ref{corolario:toro} and Theorem~\ref{Thm:Dadok}\ref{Thm:Dadok(i)}, if $\dim V_1 > 2r$ then the slice representation $(T^r, V)$ is not polar. Since $V_0$ is given by elements $X \in V$ such that $\mathrm{Ad}_g(X) = gXg^{-1} = X$ for all $g \in T^r$, it is immediate that $V_0 \subset \mathrm{Fix}(T^r,\mathfrak{sp}(2))$. Furthermore, since $\mathfrak{t}^r =\mathrm{Lie}(T^r) \subset (\mathfrak{h} + \mathfrak{k}) = V^\perp$ and $\mathfrak{t}^r \subset \mathrm{Fix}(T^r, \mathfrak{sp}(2))$, it follows that 
$$V_0 \subset \mathrm{Fix}(T^r,\mathfrak{sp}(2)) \cap (\mathfrak{t}^r)^\perp$$
and $\dim V_0 \leq \dim \mathrm{Fix}(T^r,\mathfrak{sp}(2)) - r$. Since $\mathfrak{h} \cap \mathfrak{k} = \mathfrak{t}^r$, we have $\dim (\mathfrak{h} + \mathfrak{k}) = \dim H + \dim K - r$ and then,
$$\dim V_1 = \dim V - \dim V_0 = 10 - (\dim H + \dim K - r) - \dim V_0 \geq 1 + 2r > 2r.$$
Consequently, the slice representation at $e$, given by $(T^r, V)$, is not polar and the $H \times K$-action on $\Sp(2)$ is not infinitesimally polar. 
\end{proof}

\begin{cor} Let $K$ be a closed connected subgroup of $\Sp(2)$ such that $\dim K \leq 4$ and $K$ does not contain a conjugate of $\SU(2)_{\mathrm{irr}}$. Then, the $\SU(2)_{\mathrm{irr}} \times K$-action on $\Sp(2)$ is infinitesimally polar if and only if it is almost free.\label{coro:su(2)-irr}
\end{cor}

\begin{proof}
Denote $H \coloneqq \SU(2)_{\mathrm{irr}}$. Almost-free actions are trivially infinitesimally polar, so suppose that the $H \times K$-action is not almost free, meaning that the isotropy subgroup at some $g \in \Sp(2)$ has positive dimension. Replacing $K$ by $gKg^{-1}$ (which changes neither the biquotient nor the hypotheses of the statement), we may assume that $(H \cap K)^0$ has positive dimension. Since $H$ is a simple Lie group of rank $1$ and $K$ does not contain a conjugate of $H$, the identity component $(H \cap K)^0$ must be a maximal torus of $H$. Conjugating $H$ and $K$ simultaneously by an element of $H$ if necessary, we may assume that $(H \cap K)^0 = T^1_H$, the diagonal maximal torus of $H$, which has $2$-dimensional fixed point set $\mathrm{Fix}(T^1_H,\mathfrak{sp}(2))$. Since $\dim H + \dim K \leq 7 = 9 - \dim \mathrm{Fix}(T^1_H,\mathfrak{sp}(2))$, Theorem~\ref{theo:intersecao-toro} implies that the $H \times K$-action is not infinitesimally polar.
\end{proof}

The following theorem provides a condition to determine when an action $H \times K \curvearrowright \Sp(2)$ with isotropy subgroup isomorphic to $\SU(2)$ is not infinitesimally polar.  

\begin{thm} 
Let $H$ and $K$ be closed connected subgroups of $\Sp(2)$ such that $L \coloneqq (H \cap K)^0 \cong \SU(2)$. If $\dim H + \dim K \leq 7-\dim \mathrm{Fix}(L, \mathfrak{sp}(2))$, then the $H \times K$-action on $\Sp(2)$ is not infinitesimally polar. 
\label{theo:intersecao-su2}
\end{thm}

\begin{proof} 
Let $V \coloneqq \mathfrak{h}^\perp \cap \mathfrak{k}^\perp$ be the normal space to the orbit of $e$ at this point, where the connected isotropy subgroup $L \cong \SU(2)$ acts via the slice representation $(L, V)$. Denote by $V_0$ and $V_1$ the subspaces of $V$ where $L$ acts trivially and non-trivially, respectively. By Dadok's classification \cite{dadok}, there are no polar representations of $\SU(2)$ without trivial summands of dimension greater than $5$. Hence, if $\dim V_1 > 5$, then by Theorem~\ref{Thm:Dadok}\ref{Thm:Dadok(i)} the representation $(L, V)$ is not polar. Notice that $V_0 \subset \mathrm{Fix}(L, \mathfrak{sp}(2))$, which implies that
\begin{equation*}
    \begin{aligned}
     \dim V_1 = \dim V - \dim V_0 &= 10 - (\dim H + \dim K - 3) - \dim V_0 \\ &\geq 13 - (\dim H + \dim K + \dim \mathrm{Fix}(L, \mathfrak{sp}(2))) \\ 
     &\geq 6,   
    \end{aligned}
\end{equation*}
proving that the slice representation $(L, V)$ is not polar. Consequently, the $H \times K$-action on $\Sp(2)$ is not infinitesimally polar. 
\end{proof}

\subsection{Subgroups of $\Un(2) \times \Un(2)$}\label{subsec:U2_subgroups}

For all $X \in M(4,\C)$, let $X^\ast$ and $X^t$ denote its conjugate transpose and transpose, respectively.  Under the identification of $\Sp(2)$ as a subgroup of $\SU(4)$, its Lie algebra $\mathfrak{sp}(2)$ decomposes into 
\[\mathfrak{sp}(2) = \left\{ \begin{pmatrix} 
A & B \\ 
-\thickbar{B} & \thickbar{A}
\end{pmatrix}: A^\ast = -A \ \text{and} \ B^t = B \right\} = \mathfrak{u}(2) \oplus \mathrm{Sym}(2,\C). \]
Note that this decomposition is $\Ad_{\Sp(2)}|_{\Un(2)}$-invariant, and the restriction of this representation to the subspace $\mathrm{Sym}(2,\C)$ corresponds to the action 
$$ (A,B) \in \Un(2) \times \mathrm{Sym}(2,\C) \longmapsto ABA^t \in \mathrm{Sym}(2,\C).$$

The following Lemma states that no non-trivial closed connected subgroup $H \subsetneq \Un(2)$ acts polarly on $\mathrm{Sym}(2,\C)$.

\begin{lemma} 
If $H \subsetneq \Un(2)$ is a non-trivial closed connected subgroup, then the $H$-action on $\mathrm{Sym}(2,\C)$ is not polar.  
\label{lemma:non-polar-u(2)}
\end{lemma}

\begin{proof} 
First, suppose that $H = T^2$. The $T^2$-action on $\mathrm{Sym}(2,\C)$ corresponds to 
$$\begin{pmatrix} e^{i\theta_1} & 0 \\
0 & e^{i\theta_2}
\end{pmatrix}  \begin{pmatrix} a & b \\
b & c 
\end{pmatrix} \begin{pmatrix} e^{i\theta_1} & 0 \\
0 & e^{i\theta_2}
\end{pmatrix} = \begin{pmatrix} e^{2i\theta_1}a & e^{i(\theta_1+\theta_2)}b \\
e^{i(\theta_1+\theta_2)}b & e^{2i\theta_2}c 
\end{pmatrix},$$
for $\operatorname{diag}(e^{i\theta_1}, e^{i\theta_2}, e^{-i\theta_1}, e^{-i\theta_2}) \in T^2$ and $\left(\begin{smallmatrix} a & b \\ b & c \end{smallmatrix}\right) \in \mathrm{Sym}(2,\C)$.
Hence, it turns out that $\mathrm{Sym}(2,\C)$ decomposes into weight spaces as $\C_{(2,0)} \oplus \C_{(1,1)} \oplus \C_{(0,2)}$. The cohomogeneity of the $T^2$-action on $\mathrm{Sym}(2,\C)$ is 4, whereas the cohomogeneity on each summand $\C$ is 1. Consequently, by Theorem~\ref{Thm:Dadok}\ref{Thm:Dadok(i)}, the action is not polar. Now, if $H = T^1_{m,n}$, the subspace $\mathrm{Sym}(2,\C)$ decomposes into $T^1_{m,n}$-invariant subspaces as
$$\mathrm{Sym}(2,\C) = \C_{2m} \oplus \C_{m+n} \oplus \C_{2n}.$$
Note that the cohomogeneity of this action on $\mathrm{Sym}(2,\C)$ is 5. If $m,n$ are non-zero and $m + n \neq 0$, it follows that $T^1_{m,n}$ acts non-trivially on each summand and, by Corollary \ref{corolario:toro}, this action is not polar. If either $m=-n$, or $m \cdot n = 0$, then $T^1_{m,n}$ acts trivially on one factor $\C$. In this case, the sum of the cohomogeneities on the individual factors is 4, meaning the action is still not polar. Finally, let $H = \SU(2)_{\Delta}$. Given $A \in \SU(2)$ and $B \in \mathrm{Sym}(2,\C)$, the expression $ABA^t$ corresponds to the diagonal $\SO(3)$-action on $\R^6 \cong \R^3 \oplus \R^3$. In fact, note that $\{\pm I \}$ is the ineffective kernel of the action and $\SU(2)/\{ \pm I \} \cong \SO(3)$ is its effectivization on the space $\mathrm{Sym}(2,\C) \cong \R^3 \oplus i \R^3$. This action is clearly non-polar by Theorem \ref{Thm:Dadok}\ref{Thm:Dadok(i)}, since its cohomogeneity on $\R^3 \oplus \R^3$ is $3 \neq 1 + 1$.  
\end{proof}

Let $H \subsetneq \Un(2)$ be a non-trivial closed connected subgroup. It follows that the $H \times \Un(2)$ and $H \times H$-actions on $\Sp(2)$ have isotropy subgroup at $e$ given by $H$ itself, and the normal space to the orbit of $e$ at $e$ contains the invariant subspace $\mathfrak{u}(2)^\perp = \mathrm{Sym}(2,\C)$. The corresponding slice representation is given by the restriction of $\Ad_{\Sp(2)}$ to $H$ on that normal space. By Lemma \ref{lemma:non-polar-u(2)}, the restriction of these slice representations to $\mathrm{Sym}(2,\C)$ is not polar, and owing to Theorem~\ref{Thm:Dadok}\ref{Thm:Dadok(i)}, the slice representation is likewise non-polar. Moreover, if $H,K \subsetneq \Un(2)$ are such that $H \cap K$ has positive dimension, it follows that the connected component of their intersection is conjugate to one of the $\Un(2)$ subgroups listed above. Hence, the $H \times K$-action on $\Sp(2)$ has a non-polar slice representation at $e$. In other words, the actions of the subgroups $H \times \Un(2)$, $H \times H$, and $H \times K$ with $\dim(H \cap K) > 0$ (in particular, whenever $H = T^2$) on $\Sp(2)$ are not infinitesimally polar. However, two remaining cases still need to be analyzed. For those, we proceed as follows:
\medskip

\begin{enumerate}[label=\textnormal{(\roman*)}, itemsep=0.2cm]
    \item Let $H= \SU(2)_{\Delta}$ and $K=T^1_{m,n}$. Since non-zero elements of $\mathfrak{h} = \mathrm{Lie}(H)$ and $\mathfrak{k} = \mathrm{Lie}(K)$ have eigenvalues proportional to $(1,-1,-1,1)$ and $(m,n,-m,-n)$, respectively, by Lemma \ref{lemma:almostfree}, the $H \times K$-action is almost free if and only if $\vert m \vert \neq \vert n \vert$. In this case, the biquotient $\Sp(2) \biq (\SU(2)_\Delta \times T^{1}_{m,n})$ is a 6-dimensional Riemannian orbifold. Otherwise, if $\vert m \vert = \vert n \vert$, up to conjugation, the $H \times K$-action has $T^1_{1,-1}$ as the connected component $(H \cap K)^0$. By Lemma \ref{lemma:non-polar-u(2)}, the $H \times K$-action has a non-polar slice representation at $e$ and, consequently, it is not infinitesimally polar.

    \item Let $H = T^1_{m,n}$ and $K = T^1_{a,b}$ be one-dimensional tori. There are two possibilities for the $H \times K$-action on $\Sp(2)$: either the action is almost free or, up to conjugating $K$ by some $g \in \Sp(2)$, we have $gKg^{-1} = H$, so the isotropy subgroup at $e$ is $(H \times gKg^{-1})_e = \Delta H \cong H$. In the latter case, Lemma \ref{lemma:non-polar-u(2)} implies that this action is not infinitesimally polar. In the first case, by Lemma \ref{lemma:almostfree}, we conclude that the $H \times K$-action on $\Sp(2)$ is almost free if and only if $\{ \vert a \vert, \vert b \vert \} \neq \{ \vert m \vert, \vert n \vert \}$. In this case, the biquotient $\Sp(2)\biq (T^1_{m,n} \times T^1_{a,b})$ is an $8$-dimensional Riemannian orbifold.
\end{enumerate}

\subsection{Infinitesimally polar actions on compact rank-one symmetric spaces}\label{subsec:sym-spaces}

In this subsection, we analyze product actions $H \times K \curvearrowright \Sp(2)$ where $K \in \{\SU(2) \times \Un(1), \, \SU(2) \times \SU(2), \, \SU(2)_{\low}\}$, so that $\Sp(2)/K$ is diffeomorphic to one of the compact rank-one symmetric spaces $\CP^3$, $S^4$, or $S^7$. Note that although the normal homogeneous metrics induced from $\Sp(2)$ on $\CP^3$ and $S^7$ do not coincide with the standard Fubini-Study and round metrics, the slice representation of $H_x$ at any $x \in \Sp(2)/K$ is determined up to orthogonal equivalence solely by the smooth action on $T_x(\Sp(2)/K)/T_x(H \cdot x)$. Hence, infinitesimal polarity is independent of the choice of $H$-invariant Riemannian metric on $\Sp(2)/K$, allowing us to invoke the classifications on $\CP^n$ \cite{gorodskikollross} and on spheres \cite{gorodskilyt} equipped with their standard symmetric metrics.

For the complex projective space $\CP^n$, we view it as a quotient of the unit sphere under the Hopf fibration $S^1 \hookrightarrow S^{2n+1} \to \CP^n$. Let $G$ be a compact Lie group acting by isometries on $\CP^n$ endowed with the Fubini-Study metric. Via the Hopf fibration, this action lifts to an isometric action of $\tilde{G} \coloneqq G \cdot \Un(1)$ on $\C^{n+1}$. As shown by Podest\`a and Thorbergsson \cite[Theorem~3.1]{podesta-thorbergsson}, this lifted action is polar if and only if the original $G$-action on $\CP^n$ is polar. 

\subsubsection{$H$-actions on the space $\CP^3$}\label{subsec:cp^3} The classification of infinitesimally polar actions on complex projective spaces was established by Gorodski and Kollross in \cite[Theorem~1.1]{gorodskikollross}. According to this result, an isometric action on $\CP^n$ is infinitesimally polar if and only if it is polar or orbit-equivalent to an action induced by a linear representation from an explicit list. In particular, identifying $\CP^3 \cong \Sp(2)/(\SU(2) \times \Un(1))$, the subgroups $H \subset \Sp(2)$ whose actions on $\CP^3$ are infinitesimally polar are the following:
\begin{itemize}[itemsep=0.2cm]
    \item The action of $\SU(2) \times \SU(2)$ on $\CP^3$ is polar of cohomogeneity one. Its lift
    $(\SU(2)\times\SU(2))\cdot\Un(1)$ on $\C^2\oplus\C^2$ has cohomogeneity $2$ and is
    orbit equivalent to $\Un(2)\times\Un(2)$, the isotropy representation of
    $\CP^2\times\CP^2$.

    \item The subgroup $\SU(2) \times \Un(1)$ acts polarly with cohomogeneity two on $\CP^3$. Its lift
    $\SU(2)\times\Un(1)^2$ on $\C^2\oplus\C\oplus\C$ has cohomogeneity $3$ and is orbit
    equivalent to the isotropy representation of $\CP^2\times\CP^1\times\CP^1$.

    \item The action of $\Un(2)$ is infinitesimally polar, non-polar, of cohomogeneity two. Since $\thickbar{\C}^2\cong\C^2$ as $\SU(2)$-modules, its lift is orbit equivalent to $T^2 \times \SU(2)$ on $\C^2\oplus\C^2$, which has cohomogeneity $3$ and appears in \cite[Table~1]{gorodskikollross}. It is not polar by Theorem~\ref{Thm:Dadok}\ref{Thm:Dadok(i)}, since $3\neq 1+1$.

    \item The action of $\SU(2)_\Delta$ on $\CP^3$ is infinitesimally polar, non-polar, and of cohomogeneity three. Its lift is orbit equivalent to $\Un(2)$ on $\C^2\oplus\C^2$, which has cohomogeneity $4$, also appears in \cite[Table~1]{gorodskikollross}, and is not polar by Theorem~\ref{Thm:Dadok}\ref{Thm:Dadok(i)}.
\end{itemize}

  \subsubsection{$H$-actions on the sphere $S^4$}\label{subsec:s^4}

   Next, we analyze the actions of subgroups $H$ on the sphere $S^4 \cong \Sp(2)/(\SU(2) \times \SU(2))$. Let $K$ denote $\SU(2) \times \SU(2)$. If either $H= \Un(2)$ or $H= \SU(2) \times \SU(2)$, the action $H \curvearrowright S^4$ corresponds to a cohomogeneity one Hermann action, and hence hyperpolar. Note that the case $H = \SU(2) \times \Un(1)$ was already addressed in Subsection~\ref{subsec:cp^3} up to permuting the factors. According to Gorodski and Lytchak \cite[Theorem~1.4]{gorodskilyt}, we analyze the remaining closed connected subgroups $H \subset \Sp(2)$ as follows:

   \begin{enumerate}[label=\textnormal{(\roman*)}, itemsep=0.2cm]
     
    \item\label{subsec-item4.2.2(i)} Let $H = \SU(2)_{\mathrm{low}}$ act on $S^4$. First, observe that $S^4$ naturally embeds in $\R \oplus \mathbb{H}$. Under this identification, the $H$-action on $S^4$ fixes $\R$ and acts by left multiplication on $\mathbb{H}$. Since $\Sp(1)$ acts isometrically and transitively on the spheres in $\mathbb{H}$, a regular orbit through $(x_0, q) \in S^4 \subset \R \oplus \mathbb{H}$ is the entire $3$-sphere of radius $\|q\|$ in the space $\{x_0\} \times \mathbb{H}$, so the action is polar of cohomogeneity one.

    \item Consider $H= \SU(2)_{\Delta}$. We claim that the action $H \curvearrowright S^4$ is polar. Indeed, this action corresponds to the restriction of the polar action of $\SO(3)$ on $\R^5 \cong \R^2 \oplus \R^3$ (which acts trivially on the $\R^2$ factor and transitively on the spheres in $\R^3$) to the sphere $S^4$. The principal orbits are exactly the 2-dimensional spheres in $\R^3$. Hence, the action on $S^4$ has cohomogeneity 2 and the biquotient $\Sp(2) \biq (\SU(2)_{\Delta} \times (\SU(2) \times \SU(2)))$ is a 2-dimensional Riemannian orbifold. 

    \item Let $H=T^2$. The $H$-action on $S^4$ is polar. In fact, it corresponds to the restriction of the cohomogeneity 3 polar representation of $H$ on $\R^5$, which decomposes into irreducible subspaces $\R^5 \cong \C \oplus \C \oplus \R$, to the sphere $S^4$. Hence, $H$ acts on $S^4$ with cohomogeneity $2$ and the biquotient $\Sp(2) \biq (T^2 \times (\SU(2) \times \SU(2)))$ is a 2-dimensional Riemannian orbifold. 

    \item In the case $H = \SU(2)_{\irr}$, the $H$-action on $S^4$ is polar and of cohomogeneity one. Indeed, since $K$ does not contain a conjugate of $H$, the $H$-action on $S^4$ has no fixed points on this sphere, showing that the representation of $H$ on $\R^5$ has no non-trivial fixed vectors. Hence, $\R^5$ is an irreducible $H$-representation. In other words, the $H$-action on $S^4$ is orbit equivalent to the restriction to $S^4$ of the $5$-dimensional irreducible representation of $\SO(3)$, which has cohomogeneity 2. Hence, $H$ acts with cohomogeneity one on $S^4$. 

    \item Let $H = T^1_{m,n}$. Denoting by $o = eK$ the base point of $S^4 = \Sp(2)/K$, it turns out that the isotropy representation of $H$ on $\mathfrak{p} \cong T_o(\Sp(2)/K) \cong \C^2$ has irreducible subspaces $\C_{m-n}$ and $\C_{m+n}$. If $m = \pm n$, then $H$ fixes one of the factors $\C$ and acts on the other factor $\C$, proving that the ambient space $\R^5$ of $S^4$ decomposes into $H$-invariant subspaces $\R^5 \cong \R^2 \oplus \R^3$, where $H$ fixes $\R^3$ and acts non-trivially on $\R^2$. This representation is polar and its restriction to $S^4$ is likewise polar and of cohomogeneity $3$. If $m \neq \pm n$, the representation of $H$ on $\R^5$ is a non-polar representation of cohomogeneity 4 and, by \cite[Theorem 1.4]{gorodskilyt}, the $H$-action on $S^4$ is not infinitesimally polar.
\end{enumerate}

\subsubsection{$H$-actions on the sphere $S^7$}\label{subsection:s7}
The last case to be explored in this subsection is of the quotient $\Sp(2)/\SU(2)_{\low}$. Clearly, $\SU(2)_{\low}$ is not a symmetric subgroup of $\Sp(2)$. Nevertheless, $\Sp(2)/\SU(2)_{\low}$ is diffeomorphic to the sphere $S^7 \subset \C^4$. By the classification of Gorodski and Lytchak (see \cite{gorodskilyt}), besides $H = \SU(2) \times \SU(2)$ (covered in Subsection~\ref{subsec:s^4}\ref{subsec-item4.2.2(i)}), and noting that $H = \SU(2) \times \Un(1)$ was excluded in Subsection~\ref{subsec:cp^3} while $T^2$ has isotropy $T^1$ at the basepoint $eK$ with non-polar slice representation $(T^1, \R^2 \oplus \C_1 \oplus \C_{-1})$ by Corollary~\ref{corolario:toro}, the closed connected subgroups $H$ of $\Sp(2)$ that act infinitesimally polarly on $S^7$ are the following:

 \begin{enumerate}[label=\textnormal{(\roman*)}, itemsep=0.2cm]

    \item Let $H = \Un(2)$. Owing to \cite[Theorem~1.3]{gorodskilyt}, the representation of $H$ on $\C^2 \oplus \C^2$ has infinitesimally polar restriction to $S^7$ and is not polar by Theorem~\ref{Thm:Dadok}\ref{Thm:Dadok(i)}, since the cohomogeneity on $\C^2 \oplus \C^2$ is $4 \neq 1+1$. In other words, the $H$-action on $S^7$ is infinitesimally polar and the biquotient $\Sp(2) \biq (\Un(2) \times \SU(2)_{\low})$ is a 3-dimensional Riemannian orbifold.

    \item The action of $H = \SU(2)_\Delta$ on $\Sp(2)/\SU(2)_{\low}$ is free. Indeed, if $h \in H$ fixes the point $xK \in \Sp(2)/K$, then $h \in xKx^{-1}$, where $K = \SU(2)_{\low}$. The eigenvalues of $h$ on $\C^4$ are $z,z,\thickbar{z},\thickbar{z}$ for some $z \in \Un(1)$, while those of any conjugate of an element of $\SU(2)_{\low}$ are $w,\thickbar{w},1,1$. Hence $z = 1$ and $h = e$. Consequently, the biquotient $\Sp(2)\biq (\SU(2)_\Delta \times \SU(2)_{\low})$ is a 4-dimensional Riemannian manifold. More specifically, the biquotient is a quaternionic projective space $\mathbb{H}P^1$. 

    \item Let $H = \SU(2)_{\low}$. The representation of $H$ on $\C^4 \cong \R^4 \oplus \C^2$ is trivial on $\R^4$ and acts with cohomogeneity one on $\C^2$, implying that the $H$-representation on $\C^4$ is polar. By \cite[Theorem 1.4]{gorodskilyt}, the biquotient $\Sp(2) \biq (\SU(2)_{\low} \times \SU(2)_{\low})$ is a 4-dimensional Riemannian orbifold.

    \item\label{item-action-on-s7-su(2)irr} Let $H = \SU(2)_{\irr}$. A maximal Abelian subalgebra of $\mathfrak{su}(2)_{\low}$ has weights $(0,1,0,-1)$ whereas a maximal Abelian subalgebra of $\mathfrak{h}$ has weights $(3,-1,-3,1)$. Owing to Lemma \ref{lemma:almostfree}, the $H$-action is almost free. Consequently, the biquotient $\Sp(2) \biq (\SU(2)_{\mathrm{irr}} \times \SU(2)_{\mathrm{low}})$ is a weighted quaternionic projective space $\HP^1[2]$.

    \item Let $H = T^1_{m,n}$. The space $\C^4$ decomposes into $H$-invariant subspaces as $\C^4 \cong 2\C_m \oplus 2\C_n$. If $m \cdot n \neq 0$, then the representation of $H$ on $\C^4$ is the sum of non-trivial irreducible representations (so $H$ acts almost freely on $S^7$) and, by \cite[Theorem~1.4]{gorodskilyt}, the biquotient $\Sp(2) \biq (T^1_{m,n} \times \SU(2)_{\low})$ is a $6$-dimensional Riemannian orbifold. Otherwise, suppose that $m \cdot n = 0$. In this case, by Corollary~\ref{corolario:toro}, the representation of $H$ on $\C^4$ is non-polar, has cohomogeneity $7$, and is the sum of trivial and non-trivial irreducible representations, showing by \cite[Theorem~1.4]{gorodskilyt} that the $H$-action on $S^7$ is not infinitesimally polar.
\end{enumerate}

\subsection{$H$-actions on the Berger space $B^7$} \label{section:berger-space}

In this subsection, we study actions on the 7-dimensional Berger space, beginning with an en passant review of its properties. The Berger space $B^7$ is an exceptional space in Riemannian geometry defined as the homogeneous space $B^7 \coloneqq \Sp(2)/\SU(2)_{\irr}$. It is worth mentioning that the Berger space $B^7$ is topologically different from the sphere $S^7$ and admits strictly positive curvature \cite{berger}. In order to determine infinitesimally polar actions on $B^7$, we denote $K = \SU(2)_{\irr}$ and proceed by analyzing case-by-case the remaining closed connected subgroups $H \subset \Sp(2)$ as follows.

\begin{enumerate}[label=\textnormal{(\roman*)}, itemsep=0.2cm]
    \item Let $H$ be either $\Un(2)$ or $T^2$. For both cases, at the identity, the connected component of the isotropy subgroup is isomorphic to $(H \cap K)^0 = T^1_{3,-1}$. From Corollary \ref{coro:su(2)-irr}, in both cases we conclude that the $H \times K$-action on $\Sp(2)$ is not infinitesimally polar. 
   
    \item\label{subsec:berger-o6} Let $H = T^1_{m,n}$. By Corollary \ref{coro:su(2)-irr}, the action of $H \times K$ on $\Sp(2)$ is infinitesimally polar if and only if it is almost free. Furthermore, it is clear that if $H$ and $T^1_{3,-1}$ are not conjugate, then the action is almost free. More precisely, the $H \times K$-action is almost free if and only if $\{ \vert m \vert, \vert n \vert \} \neq \{1,3\}$. In this case, the corresponding biquotient $\Sp(2)\biq (T^1_{m,n} \times \SU(2)_\mathrm{irr})$ is a $6$-dimensional Riemannian orbifold.   
    
    \item If $H = K = \SU(2)_{\irr}$, then $H \cap K = H \cong \SU(2)$ and we claim that $\mathrm{Fix}(H, \mathfrak{sp}(2)) = \{ 0 \}$. Indeed, any $X \in \mathrm{Fix}(H, \mathfrak{sp}(2))$ commutes with $H$ and thus defines an $H$-equivariant endomorphism of $\C^4$. Since $H$ acts irreducibly on $\C^4$, Schur's Lemma implies that $X = \lambda \id$ for some $\lambda \in \C$, and the condition $\mathrm{tr}(X) = 0$ yields $X = 0$. Therefore, $\mathrm{Fix}(H, \mathfrak{sp}(2)) = \{0 \}$ and, by Theorem \ref{theo:intersecao-su2}, the $H \times K$-action on $\Sp(2)$ is not infinitesimally polar.  
    
    \item\label{subsec:berger-q4} Let $H = \SU(2)_{\Delta}$. A maximal Abelian subalgebra of $\mathfrak{su}(2)_{\irr}$ has weights $(3,-1,-3,1)$ whereas a maximal Abelian subalgebra of $\mathfrak{su}(2)_{\Delta}$ has weights $(1,-1,-1,1)$. Consequently, a non-zero element in $\mathfrak{su}(2)_{\irr}$ cannot be conjugate to an element in $\mathfrak{su}(2)_{\Delta}$. By Lemma \ref{lemma:almostfree}, the $H \times K$-action on $\Sp(2)$ is almost free and the biquotient $\Sp(2) \biq (\SU(2)_{\Delta} \times \SU(2)_{\irr})$ is a 4-dimensional Riemannian orbifold. 
    \end{enumerate}
    
\subsection{Diagonal Actions} \label{seccao:diagonal}

Let $\Delta \Sp(2)$ be the diagonal subgroup of $\Sp(2) \times \Sp(2)$. The $\Delta \Sp(2)$-action on $\Sp(2)$ is hyperpolar and corresponds to the action of conjugation on $\Sp(2)$. Note that $e \in \Sp(2)$ is a fixed point of this action and the slice representation at this point is the adjoint representation $\Ad: \Sp(2) \to \GL(\mathfrak{sp}(2))$. This representation is known to be irreducible and of cohomogeneity $2$. Hence, its principal orbits have dimension $8$. 

In order to determine infinitesimally polar actions of diagonal subgroups on $\Sp(2)$, let $H$ be a closed connected subgroup of $\Sp(2)$ and denote by $\Delta H$ the diagonal subgroup of $H \times H$. With the following result, which can be found in \cite[Theorem 6]{kollrosspodesta}, we prove that no non-trivial proper connected subgroup of $\Sp(2)$ acts polarly on $\mathfrak{sp}(2)$.  

\begin{lemma} Let $\rho: G \to \On(V)$ be a faithful irreducible representation of a compact connected Lie group $G$ of cohomogeneity at least $2$. Assume the restriction of $\rho$ to a non-trivial closed connected subgroup $H$ is polar. Then the $G$- and $H$-actions on $V$ are orbit equivalent.
\label{lemma:non-polar-restriction}
\end{lemma}

For any non-trivial closed connected subgroup $H \subset \Sp(2)$, the identity $e \in \Sp(2)$ is fixed by $\Delta H \coloneqq \{(h,h) : h \in H\}$, and the slice representation at $e$ is the restriction of $\Ad_{\Sp(2)}$ to $H$ on $\mathfrak{sp}(2)$. Passing to the faithful irreducible representation of $\Sp(2)/\{\pm I_4\}$ on $\mathfrak{sp}(2)$ and the image $\thickbar{H} \coloneqq H/(H \cap \{\pm I_4\})$ (which is non-trivial and satisfies $\dim \thickbar{H} = \dim H$ since $H$ is connected), Lemma~\ref{lemma:non-polar-restriction} implies that this slice representation is polar if and only if it is orbit equivalent to $\Ad_{\Sp(2)}$, which requires $\dim H \geq 8$. By Figure~\ref{Fig:treeSp2}, this holds if and only if $H = \Sp(2)$.

\subsection{Asymmetric Subgroups}\label{subsec:asymmetric} In this subsection, we classify the asymmetric subgroups of $\Sp(2) \times \Sp(2)$, that is, subgroups that are neither a twisted diagonal nor a product, and determine which among them act infinitesimally polarly on $\Sp(2)$. To systematically construct and analyze these subgroups, we employ \textit{Goursat's Lemma}, which we presented in Section \ref{Sec:subgroupsSp2}. 

Once these asymmetric subgroups  $U \subset \Sp(2) \times \Sp(2)$ are constructed, our primary tool for testing infinitesimal polarity is the analysis of their slice representation via their weight systems. Let $U_g$ be the isotropy subgroup at $g \in \Sp(2)$ and let $T$ be its maximal torus. By translating the tangent space $T_g \Sp(2)$ back to the identity, it naturally identifies with the Lie algebra $\mathfrak{sp}(2)$, and under this identification, the action of $U_g$ corresponds to the standard adjoint action. Consequently, the complete weight system of $T$ acting on $T_g \Sp(2)$ consists exactly of the roots of $\mathfrak{sp}(2)$ restricted to $T$. Because $T_g \Sp(2)$ orthogonally decomposes into the tangent space to the orbit and the normal space, this total weight system decomposes into two disjoint multisets. Therefore, the weights of the slice representation are simply obtained by taking the restricted roots of $\mathfrak{sp}(2)$ and removing the weights associated with the tangent space to the orbit. Note that any isomorphism of $\mathfrak{su}(2) \oplus \mathfrak{su}(2)$ is realized by conjugation in $\Sp(2)$ (as swapping the two factors is inner in $\Sp(2)$), so the case $\mathfrak{h}_1 = \mathfrak{h}_2 = \mathfrak{su}(2)\oplus\mathfrak{su}(2)$ with $\mathfrak{n}_1 = \mathfrak{n}_2 = \{0\}$ reduces to the diagonal action treated in Subsection~\ref{seccao:diagonal}. Throughout this subsection and in Table~\ref{table:almost_free}, given isomorphic closed connected subgroups $H_1, H_2 \subset \Sp(2)$ and a Lie algebra isomorphism $\theta \colon \mathfrak{h}_1 \to \mathfrak{h}_2$, we denote by $\Delta^{\theta}(H_1, H_2) \coloneqq \langle \exp\{(X, \theta(X)) : X \in \mathfrak{h}_1\} \rangle \subset H_1 \times H_2$ the corresponding asymmetric diagonal subgroup (abbreviated as $\Delta^{\theta} H$ when $H_1 = H_2 = H$).

We proceed by analyzing the asymmetric closed connected subgroups of $\Sp(2) \times \Sp(2)$ case-by-case as follows:

       \subsubsection{$\mathfrak{h}_1 = \mathfrak{u}(2)$}\label{subsec:u2_goursat}
 
       In this case, the proper non-trivial ideals of $\mathfrak{h}_1$ are $\mathfrak{n}_1 = \mathfrak{su}(2)_\Delta$ (for which $\mathfrak{h}_1/\mathfrak{n}_1 \cong \mathfrak{u}(1)$) and $\mathfrak{n}_1 = \mathfrak{u}(1)$ (for which $\mathfrak{h}_1/\mathfrak{n}_1 \cong \mathfrak{su}(2)_\Delta$), while $\mathfrak{n}_1 = \{ 0 \}$ is clearly the trivial ideal of $\mathfrak{h}_1$. Moreover, in the following, we use that every $X \in \mathfrak{u}(2)$ is written as 
       \begin{align}\label{Xelementofu(2)}
        \begin{pmatrix} 
        x_{11} + i\lambda & x_{12} & 0 & 0 \\ -\thickbar{x}_{12} & -x_{11} + i\lambda & 0 & 0 \\ 0 & 0 & \thickbar{x}_{11} - i\lambda & \thickbar{x}_{12} \\ 0 & 0 & -x_{12} & - \thickbar{x}_{11} - i\lambda \end{pmatrix}
        \end{align}
        where $x_{11} \in i \R$, $x_{12} \in \C$ and $\lambda \in \R$.

        \begin{enumerate}[label=\textnormal{(\roman*)}, itemsep=0.2cm]

            \item $\mathfrak{n}_1 = \{0\}, \mathfrak{h}_2 = \mathfrak{u}(2), \mathfrak{n}_2 = \{0\}$. Notice that, in this case, every isomorphism $\theta_c : \mathfrak{h}_1/ \mathfrak{n}_1 \to \mathfrak{h}_2 /\mathfrak{n}_2$ preserves the factors of $\mathfrak{u}(2) = \mathfrak{su}(2)_\Delta \oplus \mathfrak{u}(1)$, where $\theta_c$ restricts to $\mathfrak{su}(2)_\Delta$ as the identity, up to conjugation, and restricts to $\mathfrak{u}(1)$ as a multiplication by a non-zero rational number $c$. Hence, for all $c \in \Q^\ast$, we have a $4$-dimensional subgroup $U_c$ given by 
            $$U_c \coloneqq \{ (\diag(e^{i\lambda}A, e^{-i\lambda}\thickbar{A}), \diag(e^{ic\lambda}A, e^{-ic\lambda}\thickbar{A})): A \in \SU(2), \lambda \in \R \}.$$
            Note that, if $c = \pm 1$, the subgroup $U_c$ is conjugate to the diagonal $\Delta \Un(2)$, treated in Subsection \ref{seccao:diagonal}. If $c \neq \pm 1$, it can be easily checked that the connected component of the isotropy subgroup at $e$ is such that $(U_c)^{0}_e = \Delta \SU(2)_\Delta$ and $T_e(U_c \cdot e) = \mathfrak{u}(1) \subset \mathfrak{u}(2)$, showing that the normal space $\nu_e(U_c \cdot e)$ contains the space $\mathrm{Sym}(2,\C)$. By Lemma \ref{lemma:non-polar-u(2)}, the slice representation at $e$ is not polar. 

            \item $\mathfrak{n}_1 = \mathfrak{u}(1), \mathfrak{h}_2 = \mathfrak{u}(2), \mathfrak{n}_2 = \mathfrak{u}(1)$. The subgroup $U$ that arises from this quintuple is given by
            $$U = \{ (\diag(e^{it}A, e^{-it}\thickbar{A}), \diag(e^{is}A, e^{-is}\thickbar{A})): A \in \SU(2), s,t \in \R \},$$
            whose action on $\Sp(2)$ has isotropy subgroup at $e$ given by $\Delta \Un(2) \subset \Un(2) \times \Un(2)$, showing that the tangent space $T_e(U \cdot e)$ is isomorphic to $\mathfrak{u}(1) \subset \mathfrak{u}(2)$. Hence, the normal space $\nu_e(U \cdot e)$ can be identified with $\nu_e(U \cdot e ) \cong \mathfrak{su}(2)_\Delta \oplus \Sym(2,\C)$, yielding a non-polar slice representation at $e$. In fact, the cohomogeneity of this representation on the whole space is $5$, whereas the cohomogeneities on the invariant subspaces $\mathfrak{su}(2)_\Delta$ and $\mathrm{Sym}(2,\C)$ are $1$ and $2$, respectively. Since $5 \neq 1+2$, the slice representation is not polar by Theorem~\ref{Thm:Dadok}\ref{Thm:Dadok(i)}. 

            \item \label{Subsec451item_u2_su2u1} $\mathfrak{n}_1 = \{0\}, \mathfrak{h}_2 = \mathfrak{su}(2) \oplus \mathfrak{u}(1)$ and $\mathfrak{n}_2 = \{0\}$. First, note that for this choice of $\mathfrak{h}_i$ and $\mathfrak{n}_i$, $i = 1,2$, there are infinitely many Lie algebra isomorphisms $\theta_c: \mathfrak{h}_1 \to \mathfrak{h}_2$ parametrized by constants $c \in \Q^\ast$ giving rise to a $4$-dimensional subgroup $U_c \coloneqq \Delta^{\theta_c}(\Un(2), \SU(2) \times \Un(1))$ given by
            \[
            U_c = \left\{ \left( \diag(e^{i\lambda}A, e^{-i\lambda}\thickbar{A}), \, \begin{pmatrix} \alpha & 0 & \beta & 0 \\ 0 & e^{ic\lambda} & 0 & 0 \\ -\thickbar{\beta} & 0 & \thickbar{\alpha} & 0 \\ 0 & 0 & 0 & e^{-ic\lambda} \end{pmatrix} \right) : A = \begin{pmatrix} \alpha & \beta \\ -\thickbar{\beta} & \thickbar{\alpha} \end{pmatrix} \in \SU(2), \, \lambda \in \R \right\}.
            \]
            A maximal Abelian subalgebra $\mathfrak{t}_c \subset \mathrm{Lie}(U_c)$ consists of pairs $(X,Y)$ given by
            \[
            X = \diag\big(i(\lambda + t), i(\lambda - t), -i(\lambda + t), -i(\lambda - t)\big) \qandq Y = \diag(it, ic\lambda, -it, -ic\lambda),
            \]
            with $\lambda, t \in \R$. Thus, $X$ and $Y$ are conjugate in $\Sp(2)$ if and only if they have the same eigenvalues, i.e., $\{|\lambda + t|, |\lambda - t|\} = \{|t|, |c\lambda|\}$. Clearly, this equality has no non-zero solution when $\lambda = 0$ or $t = 0$. Replacing $\lambda$ and $t$ by $|\lambda|$ and $|t|$ if necessary, we may assume $\lambda, t > 0$, so that $\lambda + t > t$. This forces $|\lambda - t| = t$ and $\lambda + t = |c|\lambda$; the first equation yields $\lambda = 2t$, and substituting into the second gives $|c| = \frac{3}{2}$. Therefore, by Lemma~\ref{lemma:almostfree}, if $|c| \neq \frac{3}{2}$, the $U_c$-action on $\Sp(2)$ is almost free and the biquotient $\Sp(2) \biq U_c$ is a $6$-dimensional Riemannian orbifold.
            
            Otherwise, if $|c| = \frac{3}{2}$, up to conjugation, we may assume $c = \frac{3}{2}$. At the point
            \[
            g = \begin{pmatrix} 0 & 1 & 0 & 0 \\ -1 & 0 & 0 & 0 \\ 0 & 0 & 0 & 1 \\ 0 & 0 & -1 & 0 \end{pmatrix} \in \Sp(2),
            \]
            the isotropy subalgebra is isomorphic to $\mathfrak{u}(1)$ and spanned by the vector $(X_0, Y_0)$, where $X_0 = \diag(3i, i, -3i, -i)$ and $Y_0 = \diag(i, 3i, -i, -3i)$. Hence, the isotropy representation of $(U_{3/2})_g^0 \cong T^1_{3,1}$ on $\mathfrak{sp}(2)$ has weights $\{0, 0, \pm 2i, \pm 2i, \pm 4i, \pm 6i\}$, whereas on the $3$-dimensional tangent space $T_g(U_{3/2} \cdot g) \cong \mathrm{Lie}(U_{3/2})/\mathrm{Lie}((U_{3/2})_g)$ it has weights $\{0, \pm 2i\}$. Consequently, $T^1_{3,1}$ acts on the $7$-dimensional normal space $\nu_g(U_{3/2} \cdot g) \cong \R \oplus \C_2 \oplus \C_4 \oplus \C_6$ with one zero weight and, by Corollary~\ref{corolario:toro} and Theorem~\ref{Thm:Dadok}\ref{Thm:Dadok(i)}, the slice representation at $g$ is not polar. 

            \item $\mathfrak{n}_1 = \mathfrak{u}(1), \mathfrak{h}_2 = \mathfrak{su}(2)\oplus\mathfrak{u}(1), \mathfrak{n}_2 = \mathfrak{u}(1)$. Using the quaternionic $2 \times 2$ matrix representation of $\Sp(2)$, the Lie algebra $\mathfrak{u}$ of the asymmetric subgroup is obtained via Goursat's Lemma by identifying the $\mathfrak{su}(2)$ factors of $\mathfrak{h}_1$ and $\mathfrak{h}_2$ through an isomorphism $\theta : \mathfrak{su}(2)_\Delta \to \mathfrak{su}(2)_{\mathrm{up}}$ mapping $\left(\begin{smallmatrix} ia & z \\ -\thickbar{z} & -ia \end{smallmatrix}\right)$ to $\diag(ia + j\thickbar{z}, 0)$. Evaluating the isotropy condition at the identity $e \in \Sp(2)$, we find that the connected component of the isotropy subgroup is the $1$-dimensional torus $ U_e^0 = \Delta T^1_{1,-1} \cong \{ \diag(e^{i\theta}, e^{-i\theta}) : \theta \in \R \}$ and the tangent space to the orbit at $e$ is given by 
            $$T_e(U \cdot e) = \left\{ \begin{pmatrix} ia - j\thickbar{z} & z \\ -\thickbar{z} & ib \end{pmatrix} : a, b \in \R, \, z \in \C \right\}.$$
            Taking the orthogonal complement inside $\mathfrak{sp}(2)$ with respect to the bi-invariant metric, the $6$-dimensional normal space is
            \[
            \nu_e(U \cdot e) = \left\{X = \begin{pmatrix} 2j\thickbar{x}_1 & x_1 + jx_2 \\ -\thickbar{x}_1 + jx_2 & jd \end{pmatrix} : x_1, x_2, d \in \C \right\}.
            \]
            For elements $h = \diag(e^{i\theta}, e^{-i\theta})$ and $X \in \nu_e(U \cdot e)$, calculating $h X h^{-1}$ yields
            \[
            \begin{pmatrix} e^{i\theta} & 0 \\ 0 & e^{-i\theta} \end{pmatrix} 
            \begin{pmatrix} 2j\thickbar{x}_1 & x_1 + jx_2 \\ -\thickbar{x}_1 + jx_2 & jd \end{pmatrix} 
            \begin{pmatrix} e^{-i\theta} & 0 \\ 0 & e^{i\theta} \end{pmatrix} 
            = \begin{pmatrix} e^{2i\theta}(2j\thickbar{x}_1) & e^{2i\theta}x_1 + jx_2 \\ -e^{-2i\theta}\thickbar{x}_1 + jx_2 & e^{-2i\theta}(jd) \end{pmatrix}.
            \]
            Hence, we obtain the $U_e^0$-invariant decomposition $\nu_e(U \cdot e) \cong \R^2 \oplus \C_2 \oplus \C_2$. Since the torus $U_e^0$ has dimension $r=1$ and acts linearly on the invariant subspace $\C_2 \oplus \C_2 \cong \C^2$ (which has complex dimension $s=2$) without non-trivial fixed vectors, Corollary~\ref{corolario:toro} implies that this restricted action is not polar. Therefore, by Theorem~\ref{Thm:Dadok}\ref{Thm:Dadok(i)}, the full slice representation at $e$ is not polar.

            \item $\mathfrak{n}_1 = \mathfrak{su}(2)_\Delta, \mathfrak{h}_2 = \mathfrak{u}(2), \mathfrak{n}_2 = \mathfrak{su}(2)_\Delta$. Given the Lie subalgebras $\mathfrak{h}_1, \mathfrak{n}_1, \mathfrak{h}_2, \mathfrak{n}_2$, every isomorphism $\theta_c: \mathfrak{h}_1/\mathfrak{n}_1 \to \mathfrak{h}_2/\mathfrak{n}_2$ is given by multiplication by $c = m/k \in \Q^\ast$, where $m,k$ are coprime, and gives rise to a subgroup
            \[
            U_c = \{ (\diag(e^{ikt}A, e^{-ikt}\thickbar{A}), \diag(e^{imt}B, e^{-imt}\thickbar{B})): A, B \in \SU(2), \, t \in \R \}.
            \]
            If $c \neq 1$, it follows that the connected component of the isotropy subgroup at $e$ is $(U_c)_e^0 = \Delta \SU(2)_\Delta \subset \SU(2)_\Delta \times \SU(2)_\Delta$, and the corresponding normal space $\nu_e(U_c \cdot e)$ is $\Sym(2,\C)$. If $c = 1$, at the point 
            \[
            g = \begin{pmatrix}
             0 & I_2 \\
             -I_2 & 0
            \end{pmatrix} \in \Sp(2),
            \]
            the connected component of the isotropy subgroup is $(U_1)_g^0 \cong \SU(2)_\Delta$ and the normal space $\nu_g(U_1 \cdot g)$ is likewise identified with $\Sym(2,\C)$. For both cases, Lemma~\ref{lemma:non-polar-u(2)} implies that these slice representations are not polar and the $U_c$-action on $\Sp(2)$ is not infinitesimally polar for all values of the non-zero rational $c$.

            \item $\mathfrak{n}_1 = \mathfrak{su}(2)_\Delta, \mathfrak{h}_2= \mathfrak{su}(2) \oplus \mathfrak{u}(1)$ and $\mathfrak{n}_2= \mathfrak{su}(2) \oplus \{0\}$. Since $\mathfrak{h}_1/\mathfrak{n}_1 \cong \mathfrak{h}_2/\mathfrak{n}_2 \cong \R$, the isomorphism $\theta$ from Goursat's Lemma corresponds to a multiplication by some non-zero rational constant $c \in \Q^\ast$. The subalgebra $\mathfrak{u}$, as given by \eqref{Goursatslemma:subalgu}, consists of all pairs of matrices $(X,Y) \in \mathfrak{sp}(2) \oplus \mathfrak{sp}(2)$ such that $X$ is given as in \eqref{Xelementofu(2)} and        
        \[Y = \begin{pmatrix} y_{11} & 0 & y_{13} & 0 \\ 0 & i(c\lambda) & 0 & 0 \\ -\thickbar{y}_{13} & 0 & \thickbar{y}_{11} & 0 \\ 0 & 0 & 0 & -i(c\lambda) \end{pmatrix}
        \]
        where $y_{11} \in i \R$, $y_{13} \in \C$ and $\lambda \in \R$. Consequently, $\dim \mathfrak{u}=7$. The isotropy subgroup of the $U$-action on $\Sp(2)$ at the identity has Lie algebra isomorphic to $\mathfrak{u}(1)$ and is spanned by the element $Z= \diag(z_1,z_2, -z_1,-z_2)$ where $z_1=i(2-c)$ and $z_2=ic$. The tangent space to the $U$-orbit at the identity has dimension $7-1=6$ and its normal space consists of symmetric off-diagonal block matrices of the form
        \[
        \begin{pmatrix}
            0 & x \\
            x & y
        \end{pmatrix} \quad \textnormal{with} \quad x,y \in \C .
        \]
        The slice representation corresponds to $\Un(1) \curvearrowright \R^4 \cong \C^2$ with non-zero weights $z_1 + z_2 = 2i$ and $2z_2 = 2ic$. Owing to Corollary \ref{corolario:toro}, it is not polar. Therefore, the $U$-action is not infinitesimally polar.

        \item $\mathfrak{n}_1 = \mathfrak{su}(2)_\Delta, \mathfrak{h}_2=\mathfrak{t}^2$ and $\mathfrak{n}_2=\mathfrak{t}^1_{m,n}$. Since $\gcd(m,n)=1$, it follows that $\mathfrak{h}_2/\mathfrak{n}_2 \cong \R$. The corresponding projection associates each $Y = \diag(i \phi_1, i \phi_2, -i \phi_1, - i \phi_2)$ to the real number $m \phi_2 - n \phi_1 \in \R$. The isomorphism $\theta$ is a multiplication by a non-zero rational number $c \in \Q^\ast$. The fiber product condition translates to $m \phi_2 - n \phi_1 = c \lambda$ and the subalgebra $\mathfrak{u}$ has dimension $4+2-1=5$. Since $U \subset \Un(2) \times T^2 \subset \Un(2) \times \Un(2)$ and $\dim \mathfrak{u} = 5$, we have $T_e(U \cdot e) \subset \mathfrak{u}(2)$, which implies that $1 \leq \dim U_e \leq 2$, with $U_e^0 \subset \Delta T^2 \subsetneq \Delta\Un(2)$, and $\mathrm{Sym}(2,\C) = \mathfrak{u}(2)^\perp \subset \nu_e(U \cdot e)$. By Lemma \ref{lemma:non-polar-u(2)} and Theorem \ref{Thm:Dadok}\ref{Thm:Dadok(i)}, the slice representation at $e$ is not polar, so the $U$-action is not infinitesimally polar.

        \item \label{Subsec451itemiii} $\mathfrak{n}_1 = \mathfrak{su}(2)_\Delta, \mathfrak{h}_2=\mathfrak{t}^1_{m,n}$ and $\mathfrak{n}_2=\{0\}$. The quotient for the first factor is $\mathfrak{h}_1/\mathfrak{n}_1 \cong \R$, with the projection $\pi_1(X) = \lambda \in \R$ for $X \in \mathfrak{u}(2)$ as in \eqref{Xelementofu(2)}. For the second factor, since the ideal is trivial, the quotient is $\mathfrak{h}_2/\mathfrak{n}_2 \cong \R$. The corresponding projection $\pi_2 : \mathfrak{t}^1_{m,n} \to \R$ associates each $Y = \diag(imt, int, -imt, -int)$ with the real number $\pi_2(Y) = t \in \R$. The isomorphism $\theta$ corresponds to multiplication by a non-zero rational constant $c \in \Q^\ast$. The fiber product condition $\theta(\pi_1(X)) = \pi_2(Y)$ translates simply to $c \lambda = t$. The subalgebra $\mathfrak{u}$ has dimension $4+1-1=4$ and consists of pairs $(X,Y)$ where 
        \[
        Y = \diag\left(im(c\lambda), in(c\lambda), -im(c\lambda), -in(c\lambda) \right)  
        \]
        for $\lambda \in \R$, generating the $4$-dimensional subgroup $U^c_{m,n} \subset \Un(2) \times T^1_{m,n}$ given by
        \[
        U^c_{m,n} = \left\{ \left( \diag(e^{i\lambda}A, e^{-i\lambda}\thickbar{A}), \, \diag(e^{imc\lambda}, e^{inc\lambda}, e^{-imc\lambda}, e^{-inc\lambda}) \right) : A \in \SU(2), \, \lambda \in \R \right\}.
        \]
        By Lemma \ref{lemma:almostfree}, if $c(\pm m \pm n) \neq 2$ (for all choices of signs), the isotropy subalgebra is trivial at every point of $\Sp(2)$, so the $U^c_{m,n}$-action is almost free and the biquotient $\Sp(2) \biq U^c_{m,n}$ is a $6$-dimensional Riemannian orbifold. 
        
        Otherwise, up to conjugation by the Weyl group (which changes the signs of $m$ and $n$), we may assume $m+n \neq 0$ and $c = \frac{2}{m+n}$, so it suffices to analyze the isotropy at the identity. The isotropy subalgebra at the identity consists of elements $(X,Y) \in \mathfrak{u}$ such that $X = Y$. Equating the blocks of \eqref{Xelementofu(2)} with $Y$ yields $x_{12} = 0$. Writing $x_{11} = ia$ for $a \in \R$, the diagonal entries give $a + \lambda = mc\lambda$ and $-a + \lambda = nc\lambda$, so $a = \lambda \frac{m-n}{m+n}$ and the connected isotropy subgroup at $e$ is $(U^c_{m,n})_e^0 = \Delta T^1_{m,n} \cong T^1_{m,n}$. Since $T_e(U^c_{m,n} \cdot e) = \mathfrak{su}(2)_\Delta \subset \mathfrak{u}(2)$, the normal space $\nu_e(U^c_{m,n} \cdot e)$ contains $\mathfrak{u}(2)^\perp = \mathrm{Sym}(2,\C)$, and by Lemma \ref{lemma:non-polar-u(2)} the slice representation at $e$ is not polar.

        \item $\mathfrak{n}_1 = \mathfrak{u}(1), \mathfrak{h}_2= \mathfrak{su}(2)_{\irr}$ and $\mathfrak{n}_2=\{ 0\}$. The Lie subalgebra $\mathfrak{u}$ generated by Goursat's Lemma is given by $\mathfrak{u} = \{ (X_s+X_c, \theta(X_s)): X_s \in \mathfrak{su}(2)_\Delta, X_c \in \mathfrak{u}(1) \},$
        where $\theta: \mathfrak{su}(2)_\Delta \to \mathfrak{su}(2)_{\irr}$ is an isomorphism. At the point 
        $$g = \begin{pmatrix} 1 & 0 & 0 & 0 \\ 0 & 0 & 0 & 1 \\ 0 & 0 & 1 & 0 \\ 0 & -1 & 0 & 0 \end{pmatrix} \in \Sp(2),$$
        the isotropy subgroup of the action on $\Sp(2)$ has Lie algebra isomorphic to $\mathfrak{u}(1)$ and spanned by the vector $(X_1, X_2) \in \mathfrak{u}$ such that $X_1 = \diag(3i, i, -3i, -i)$ and $X_2 = \diag (3i,-i,-3i, i)$. Since the weights of the corresponding isotropy representation on $\mathfrak{sp}(2)$ are given in the set $\{ 0, 0 , \pm 2i, \pm 2i, \pm 4i, \pm 6i \}$, it straightforwardly follows that the slice representation at $g$ is non-polar. In fact, the normal space to the orbit of $g$ at this point is $7$-dimensional and it has the weight zero with multiplicity one. That is, this normal space decomposes into $\R \oplus \C_2 \oplus \C_4 \oplus \C_6$. By Corollary~\ref{corolario:toro} and Theorem~\ref{Thm:Dadok}\ref{Thm:Dadok(i)}, the action on $\Sp(2)$ is not infinitesimally polar. 

        \item $\mathfrak{n}_1 = \mathfrak{u}(1), \mathfrak{h}_2= \mathfrak{su}(2) \oplus \mathfrak{su}(2)$ and $\mathfrak{n}_2= \mathfrak{su}(2)_{\low}$. The subgroup $U$ associated to the quintuple $\{\mathfrak{h}_1, \mathfrak{n}_1, \mathfrak{h}_2, \mathfrak{n}_2, \theta \}$ is given by 
        $$U = \left\{ \left( \diag(e^{it}A_1, e^{-it}\thickbar{A}_1) , \begin{pmatrix} \alpha_1 & 0 & \beta_1 & 0 \\ 0 & \alpha_2 & 0 & \beta_2 \\ -\thickbar{\beta}_1 & 0 & \thickbar{\alpha}_1 & 0 \\ 0 & -\thickbar{\beta}_2 & 0 & \thickbar{\alpha}_2
            
        \end{pmatrix} \right) : A_i = \begin{pmatrix}
            \alpha_i & \beta_i \\
            -\thickbar{\beta}_i & \thickbar{\alpha}_i 
        \end{pmatrix} \in \SU(2), t \in \R \right\}.$$
        At the point 
        $$g = \begin{pmatrix} 0 & 1 & 0 & 0 \\ -1 & 0 & 0 & 0 \\ 0 & 0 & 0 & 1 \\ 0 & 0 & -1 & 0 \end{pmatrix} \in \operatorname{Sp}(2),$$
        the isotropy subgroup consists of the elements $\alpha_1, \alpha_2 \in \Un(1)$ and $\beta_1 = \beta_2 = 0$ such that $\alpha_2 = \alpha_1^3$. Hence, the connected component of the isotropy subgroup $U_g$ can be identified with
        $$U_g^0 = \left\{ \left( \operatorname{diag}(e^{3it}, e^{it}, e^{-3it}, e^{-it}), \; \operatorname{diag}(e^{it}, e^{3it}, e^{-it}, e^{-3it}) \right) : t \in \mathbb{R} \right\}.$$
        It is clear that the isotropy representation at $g$ has weights $\{0, 0, \pm 2i, \pm 2i, \pm 4i, \pm 6i \}$. Since the Lie algebra $\mathfrak{u}$ of $U$ decomposes into $\mathfrak{u} \cong \mathfrak{n}_1 \oplus \mathfrak{n}_2 \oplus \gr(\theta),$
        where $\gr(\theta) \cong \mathfrak{su}(2)$ is the graph of $\theta$, it follows that zero is a weight of multiplicity two with respect to the isotropy representation restricted to the tangent space $T_g (U \cdot g)$. Consequently, the slice representation $(U_g^0, \C^2)$ has only non-zero weights and we conclude, by Corollary~\ref{corolario:toro}, that the $U$-action on $\Sp(2)$ is not infinitesimally polar.

        \item $\mathfrak{n}_1 = \mathfrak{u}(1), \mathfrak{h}_2=\mathfrak{su}(2)_{\Delta}$ and $\mathfrak{n}_2=\{0\}$. In this case, the subgroup $U$ associated to the quintuple $\{\mathfrak{h}_1, \mathfrak{n}_1, \mathfrak{h}_2, \mathfrak{n}_2, \id \}$ is given by
        $$U = \{ (\diag(e^{i\theta}A, e^{-i\theta}\thickbar{A}), \diag(A, \thickbar{A})): A \in \SU(2), \theta \in \R \}.$$
        At the identity $e \in \Sp(2)$, the isotropy subgroup $U_e$ is given by $U_e = \Delta (\SU(2)_\Delta) \subset \Sp(2) \times \Sp(2)$. Notice that the Lie algebra $\mathfrak{u}$ of $U$ is given by $\mathfrak{u} = \mathfrak{n}_1 \oplus \Delta (\mathfrak{su}(2)_\Delta),$
        and hence $T_e(U \cdot e) \cong \mathfrak{n}_1 \subset \mathfrak{u}(2)$. Hence, the normal space $\nu_e ( U \cdot e)$ contains $\mathfrak{u}(2)^\perp = \mathrm{Sym}(2,\C)$, and by Lemma \ref{lemma:non-polar-u(2)} the slice representation at $e$ is not polar. 

       \item \label{Subsec451itemvii}$\mathfrak{n}_1 = \mathfrak{u}(1), \mathfrak{h}_2=\mathfrak{su}(2)_{\low}$ and $\mathfrak{n}_2=\{0\}$. This quintuple gives rise to the subalgebra $\mathfrak{u}$ given by
        $$\mathfrak{u} = \left\{ \left( 
            \begin{pmatrix} 
            i(t+a) & z & 0 & 0 \\ 
            -\thickbar{z} & i(t-a) & 0 & 0 \\ 
            0 & 0 & -i(t+a) & \thickbar{z} \\ 
            0 & 0 & -z & -i(t-a) 
            \end{pmatrix}, \,
            \begin{pmatrix} 
            0 & 0 & 0 & 0 \\ 
            0 & ia & 0 & z \\ 
            0 & 0 & 0 & 0 \\ 
            0 & -\thickbar{z} & 0 & -ia 
            \end{pmatrix} 
            \right) : t, a \in \mathbb{R}, z \in \mathbb{C} \right\},$$
        which we claim corresponds to an almost-free action on $\Sp(2)$. In fact, setting $\lambda = \sqrt{a^2 + |z|^2}$, a matrix $X_1$ on the left factor of $\mathfrak{u}$ has eigenvalues $\{i(t+\lambda), -i(t+\lambda), i(t-\lambda), -i(t-\lambda)\}$, whereas a matrix $X_2$ on its right factor has eigenvalues $\{ 0, 0, i\lambda, -i\lambda \}$. Hence, $X_1 = \Ad_g(X_2)$ for some $g \in \Sp(2)$ only if they have the same eigenvalues. However, if they are in fact conjugate, it follows that $|t| = \lambda$, so that the eigenvalues of $X_1$ reduce to $\{0, 0, 2i\lambda, -2i\lambda\}$. Comparing both spectra yields $2\lambda = \lambda$, which implies $\lambda = 0$ and hence $t = a = z = 0$. Therefore, the isotropy subalgebra is trivial at every point of $\Sp(2)$, so the action is almost free and the biquotient $\Sp(2) \biq U$, where $U = \Un(1) \cdot \Delta^{\theta}(\SU(2)_{\Delta}, \SU(2)_{\low})$, is a $6$-dimensional Riemannian orbifold. 
        \end{enumerate}         
    
   \subsubsection{$\mathfrak{h}_1 = \mathfrak{su}(2)_{\irr}$, $\mathfrak{n}_1 = \{ 0 \}$}\label{subsec:su2irr_goursat}

   \begin{enumerate}[label=\textnormal{(\roman*)}, itemsep=0.2cm]
       \item $\mathfrak{h}_2= \mathfrak{su}(2)_\Delta$ and $\mathfrak{n}_2 = \{0\}$. Since $\SU(2)_{\irr} \times \SU(2)_\Delta$ acts almost freely on $\Sp(2)$ (see Subsection \ref{section:berger-space}\ref{subsec:berger-q4}), its subgroup $\Delta^{\theta}(\SU(2)_{\irr}, \SU(2)_{\Delta})$ also acts almost freely on $\Sp(2)$, and the biquotient $\Sp(2) \biq \Delta^{\theta}(\SU(2)_{\irr}, \SU(2)_{\Delta})$ is a $7$-dimensional Riemannian orbifold (denoted $\Sp(2)\biq \SU(2)_{\varphi_3}$ in \cite[Subsec.~5.2.2]{Andregomesthesis}).

       \item $\mathfrak{h}_2 = \mathfrak{su}(2) \oplus \mathfrak{su}(2)$ and $\mathfrak{n}_2 = \mathfrak{su}(2)_{\low}$. The corresponding $6$-dimensional subgroup $U$ of $\SU(2)_{\irr} \times (\SU(2) \times \SU(2))$ acts on $\Sp(2)$ with isotropy subgroup at the point 
       $$q = \begin{pmatrix} 0 & 1 & 0 & 0 \\ 0 & 0 & -1 & 0 \\ 0 & 0 & 0 & 1 \\ 1 & 0 & 0 & 0 \end{pmatrix} \in \Sp(2)$$
        whose connected component is, up to conjugation in $\Sp(2) \times \Sp(2)$, given by the diagonal $\Delta T^1_{3,-1}$ of a maximal torus of $\SU(2)_{\irr}$. Since $\Sp(2)$ has roots $\pm 2e_1, \pm 2e_2, \pm e_1 \pm e_2$ and the corresponding weights of $T^1_{3,-1}$ are $\lambda_1 = 3i$ and $\lambda_2 = -i$, it follows that the isotropy representation $(T^1_{3,-1}, \mathfrak{sp}(2))$ has weights $\{0,0, \pm2i,  \pm2i, \pm 4i, \pm 6i\}$. Furthermore, since the Lie algebra $\mathfrak{u} \coloneqq \mathrm{Lie}(U)$ is isomorphic to $\mathfrak{su}(2) \oplus \mathfrak{su}(2)$, it follows that the representation $(T^1_{3,-1}, T_q (U \cdot q))$ has weights $\{0, \pm 2i, \pm 6i\}$. Consequently, the weights of the slice representation $(T^1_{3,-1}, \nu_q (U \cdot q))$ are $\{0, \pm 2i, \pm 4i\}$. Hence, by Theorem~\ref{Thm:Dadok}\ref{Thm:Dadok(i)} and Corollary \ref{corolario:toro}, the slice representation at $q$ is not polar and the $U$-action on $\Sp(2)$ is not infinitesimally polar.

       \item $\mathfrak{h}_2 = \mathfrak{su}(2) \oplus \mathfrak{u}(1)$ and $\mathfrak{n}_2 = \{0\} \oplus \mathfrak{u}(1)$. The corresponding $4$-dimensional subgroup $U$ of $\SU(2)_{\irr} \times (\SU(2) \times \Un(1))$ acts on $\Sp(2)$ with the same isotropy subgroup $T^1_{3,-1}$ at $q \in \Sp(2)$ as in the previous item. Since $\mathfrak{u} \coloneqq \mathrm{Lie}(U)$ is isomorphic to $\mathfrak{su}(2) \oplus \mathfrak{u}(1)$, the representation $(T^1_{3,-1}, T_q(U \cdot q))$ has weights $\{0, \pm 2i\}$, and consequently the weights of the slice representation $(T^1_{3,-1}, \nu_q(U \cdot q))$ are $\{0, \pm 2i, \pm 4i, \pm 6i\}$. Hence, by Theorem~\ref{Thm:Dadok}\ref{Thm:Dadok(i)} and Corollary~\ref{corolario:toro}, the slice representation at $q$ is not polar and the $U$-action on $\Sp(2)$ is not infinitesimally polar.

       \item $\mathfrak{h}_2=\mathfrak{su}(2)_{\low}$ and $\mathfrak{n}_2=\{0\}$. Since $\SU(2)_{\irr} \times \SU(2)_{\low}$ acts almost freely on $\Sp(2)$ (see Subsection \ref{subsection:s7}\ref{item-action-on-s7-su(2)irr}), its subgroup $\Delta^{\theta}(\SU(2)_{\irr}, \SU(2)_{\low})$ also acts almost freely on $\Sp(2)$, and the biquotient $\Sp(2) \biq \Delta^{\theta}(\SU(2)_{\irr}, \SU(2)_{\low})$ is a $7$-dimensional Riemannian orbifold (denoted $\Sp(2)\biq \SU(2)_{\varphi_2}$ in \cite[Subsec.~5.2.2]{Andregomesthesis}).
   \end{enumerate} 

      \subsubsection{$\mathfrak{h}_1=\mathfrak{su}(2) \oplus \mathfrak{su}(2), \mathfrak{n}_1=\mathfrak{su}(2)_{\low}$}\label{subsec:su2xsu2_goursat}

      In the following cases, we identify the $\Sp(2)$ maximal subgroup $\SU(2) \times \SU(2)$ as $\Sp(1) \times \Sp(1) = \{ \diag(p,q) : p,q \in \Sp(1) \}$. Similarly, consider the identifications of the subgroups $\SU(2)$ as 
      $$\SU(2)_{\Delta} \cong \Sp(1)_{\Delta} \coloneqq \{ \diag(q,q) : q \in \Sp(1) \}$$
      and 
      $$\SU(2)_{\low} \cong \Sp(1)_{\low} = \{ \diag(1,w): w \in \Sp(1) \}.$$
       In addition, it is clear that the subgroup $\SU(2)_{\low}$ is conjugate in $\Sp(2)$ to the upper-block subgroup given by $\SU(2)_{\mathrm{up}} \cong \Sp(1)_{\mathrm{up}} \coloneqq \{ \diag(w,1) : w \in \Sp(1) \}$. In cases~(ii)--(iv) below, the tangent space to the orbit at the identity $e \in \Sp(2)$ is $T_e (U \cdot e) \cong \mathfrak{sp}(1)_{\low}$, and hence the $7$-dimensional normal space $\nu_e(U \cdot e)$ is given by
    \begin{equation}
        \begin{aligned}
            \nu_e (U \cdot e) = \left\{ \begin{pmatrix} u & -\thickbar{x} \\
    x & 0 
    \end{pmatrix}: u \in \mathrm{Im}(\mathbb{H}), x \in \mathbb{H} \right\}.
        \end{aligned}
        \label{eq:normal-space-goursat}
    \end{equation}

      \begin{enumerate}[label=\textnormal{(\roman*)}, itemsep=0.2cm]

         \item\label{Subsec453itemi} $\mathfrak{h}_2 = \mathfrak{su}(2) \oplus \mathfrak{su}(2), \mathfrak{n}_2 = \mathfrak{su}(2)_{\low}$. This choice of quintuple gives rise to the subgroup $U \coloneqq \Delta\SU(2)_{\mathrm{up}} \cdot (\SU(2)_{\low} \times \SU(2)_{\low})$ given by
         $$U = \{(\diag(p, q_1), \diag(p, q_2)) : p, q_1, q_2 \in \Sp(1)\} \cong \Sp(1)^3,$$
         via the identification $\SU(2) \cong \Sp(1)$. Observe that $N_2 = \{ (\id , \diag(1,q_2)) : q_2 \in \Sp(1) \} \subset U$ acts on $\Sp(2)$ with quotient $\Sp(2)/N_2 \cong S^7$. Furthermore, the induced action of $U/N_2 \cong \Sp(1) \times \Sp(1)$ on $S^7 \subset \mathbb{H}^2$ is given by $(p,q_1) \cdot (a,c) = (pap^{-1}, q_1cp^{-1})$, which is orbit equivalent to the restriction to $S^7$ of the cohomogeneity $3$ polar action $\SO(3) \times \SO(4) \curvearrowright \R^8 \cong \R \oplus \R^3 \oplus \R^4$, where $\R$ is the fixed component. Consequently, the $U$-action on $\Sp(2)$ is infinitesimally polar and the biquotient $\Sp(2) \biq U$ is a 2-dimensional Riemannian orbifold.
                     
          \item $\mathfrak{h}_2 = \mathfrak{su}(2) \oplus \mathfrak{u}(1), \mathfrak{n}_2 = \{0\} \oplus \mathfrak{u}(1)$. This quintuple gives rise to the subgroup 
      $$U = \{ (\diag(p,q), \diag(p,z)): p,q \in \Sp(1), z \in \Un(1) \},$$
      which acts on $\Sp(2)$ with isotropy subgroup $U_e$ at the identity $e$ given by
      $$U_e = \{ (\diag(p,z), \diag(p,z)): p \in \Sp(1), z \in \Un(1) \} \cong \Sp(1) \times \Un(1).$$
      Consequently, at $e$, the slice representation corresponds to $U_e$ acting on $\nu_e(U \cdot e) \cong \R^3 \oplus \C^2$, where the action of $U_e$ on $\C^2$ is orbit equivalent to the standard $\Un(2)$-action and $\Sp(1)$ acts on $\R^3$ via the $\Ad$-representation of $\SO(3)$. Since its cohomogeneity is $3 \neq 1+1$, by Theorem~\ref{Thm:Dadok}\ref{Thm:Dadok(i)}, the slice representation at $e$ is not polar. 

      \item $\mathfrak{h}_2= \mathfrak{su}(2)_\Delta$ and $\mathfrak{n}_2=\{0\}$. For this choice of $\mathfrak{h}_2$ and $\mathfrak{n}_2$, the associated subgroup $U$ is given by 
   $U = \{ (\diag(p,q), \diag(p,p)): p,q \in \Sp(1) \}.$
   At the identity $e \in \Sp(2)$, the isotropy subgroup is given by $U_e = \{ (\diag(p,p), \diag(p,p)): p \in \Sp(1) \} = \Delta \Sp(1)_\Delta,$
    implying that the normal space \eqref{eq:normal-space-goursat} decomposes into $U_e$-irreducible subspaces as $\nu_e(U \cdot e) \cong \R \oplus \R^3 \oplus \R^3$, where $U_e$ acts on both factors $\R^3$ as the $\SO(3)$ $\Ad$-representation and it acts trivially on $\R$, showing that the slice representation at $e$ is not polar and the $U$-action on $\Sp(2)$ is not infinitesimally polar.

    \item $\mathfrak{h}_2= \mathfrak{su}(2)_{\low}$ and $\mathfrak{n}_2=\{0\}$. Notice that $\Sp(1)_{\low}$ is conjugate in $\Sp(2)$ to $\Sp(1)_{\mathrm{up}}$ and the quintuple $\{ \mathfrak{h}_1, \mathfrak{n}_1, \mathfrak{h}_2, \mathfrak{n}_2, \theta \}$ gives rise to the group $U = \{ (\diag(p,q), \diag(p,1)): p,q \in \Sp(1) \}.$
    It can be easily checked that $U$ acts on $\Sp(2)$ with isotropy subgroup at $e$ given by 
    $$U_e = \{ (\diag(p,1), \diag(p,1)): p \in \Sp(1) \} = \Delta \Sp(1)_{\mathrm{up}},$$
    implying that the normal space \eqref{eq:normal-space-goursat} decomposes into $U_e$-irreducible subspaces as $\nu_e (U \cdot e ) \cong \R^3 \oplus \C^2$, where $U_e$ acts on $\R^3$ as the $\SO(3)$ $\Ad$-representation and it acts on $\C^2$ as the standard $\SU(2)$-representation. Consequently, the slice representation at $e$ is not polar and $U$ does not act infinitesimally polarly on $\Sp(2)$.   

      \end{enumerate}

    \subsubsection{$\mathfrak{h}_1=\mathfrak{su}(2)_{\low}, \mathfrak{n}_1=\{0\}, \mathfrak{h}_2=\mathfrak{su}(2)_\Delta, \mathfrak{n}_2=\{0\}$} \label{subsec:GM_sphere}

    By Goursat's Lemma, this case yields the subgroup $\Delta^{\theta}(\SU(2)_{\low}, \SU(2)_{\Delta})$, whose biquotient $\Sp(2) \biq \Delta^{\theta}(\SU(2)_{\low}, \SU(2)_{\Delta})$ is the Gromoll-Meyer sphere \cite{gromoll}, which is a $7$-dimensional Riemannian manifold.

    \subsubsection{$\mathfrak{h}_1=\mathfrak{u}(1) \oplus \mathfrak{su}(2)$}\label{subsec:U_mk} In the corresponding cases below, we still consider $\Sp(2)$ and its subgroups with the quaternionic expression. 

    \begin{enumerate}[label=\textnormal{(\roman*)}, itemsep=0.2cm]

        \item\label{Subsec455item_t1} $\mathfrak{n}_1 = \mathfrak{su}(2)_{\low}, \mathfrak{h}_2 = \mathfrak{t}^1_{m,n}, \mathfrak{n}_2 = \{0\}$. By Goursat's Lemma, for each $c \in \Q^\ast$, this quintuple gives rise to the $4$-dimensional subgroup
        \[
        V^c_{m,n} \coloneqq \{ (\diag(e^{ict}, p), \diag(e^{imt}, e^{int})) : t \in \R, \, p \in \Sp(1) \}.
        \]
        Since $N_1 = \{ (\diag(1,p), e) : p \in \Sp(1) \} \subset V^c_{m,n}$ acts freely on $\Sp(2)$ with quotient $\Sp(2)/N_1$ diffeomorphic to the sphere $S^7 \subset \mathbb{H}^2$, the orbit space $\Sp(2) \biq V^c_{m,n}$ coincides with that of $V^c_{m,n}/N_1 \cong S^1$ acting on $S^7$ via $(a,b) \mapsto (e^{ict} a e^{-imt}, e^{ict} b e^{-int})$. Writing $a = z_1 + z_2 j$ and $b = z_3 + z_4 j$ ($z_k \in \C$), the induced $S^1$-representation on $\C^4$ has weights $c - m$, $c + m$, $c - n$, and $c + n$. If $|c| \notin \{|m|, |n|\}$, all four weights are non-zero, so the action is almost free and the biquotient $\Sp(2) \biq V^c_{m,n}$ is a $6$-dimensional Riemannian orbifold. Otherwise, if $|c| \in \{|m|, |n|\}$, at least one weight vanishes while at least two remain non-zero (as $(c-m)+(c+m) = (c-n)+(c+n) = 2c \neq 0$). Hence, the $S^1$-representation on $\C^4$ is the sum of trivial and non-trivial summands and is non-polar by Corollary~\ref{corolario:toro} and Theorem~\ref{Thm:Dadok}\ref{Thm:Dadok(i)}, showing by \cite[Theorem~1.4]{gorodskilyt} that the $V^c_{m,n}$-action on $\Sp(2)$ is not infinitesimally polar.

        \item $\mathfrak{n}_1 = \{0\}, \mathfrak{h}_2 = \mathfrak{u}(1) \oplus \mathfrak{su}(2), \mathfrak{n}_2 = \{0\}$. In this case, the Lie algebra of the asymmetric subgroup is the graph of an isomorphism between $\mathfrak{h}_1$ and $\mathfrak{h}_2$. Assuming the isomorphism is the identity on the $\mathfrak{su}(2)$ factor and scales the $\mathfrak{u}(1)$ factor by a constant $c \in \Q^\ast$, we obtain
        $$\mathfrak{u} = \{ (\diag(i\alpha, q), \diag(ic\alpha, q)) : \alpha \in \R, q \in \mathfrak{sp}(1) \}.$$ 
        Note that for $c = \pm 1$, the subgroup is conjugate to $\Delta(\Un(1) \times \SU(2))$, treated in Subsection~\ref{seccao:diagonal}, and we now analyze the case $c \neq \pm 1$. The isotropy subalgebra at $e \in \Sp(2)$ is determined by the elements satisfying $\diag(i\alpha, q) = \diag(ic\alpha, q)$. Since $c \neq 1$, this forces $\alpha = 0$, leaving $q \in \mathfrak{sp}(1)$ free. It follows that
        \[\mathfrak{u}_e = \{ (\diag(0, q), \diag(0, q)) : q \in \mathfrak{sp}(1) \} \cong \Delta \mathfrak{sp}(1)_{\low}\]
        and the connected component of the isotropy subgroup is exactly $U_e^0 = \Delta\Sp(1)_{\low}$. The tangent space to the orbit at $e$ is given by $T_e(U \cdot e) = \{ \diag(i(1-c)\alpha, 0) : \alpha \in \R \} \cong \mathfrak{u}(1)$. Taking the orthogonal complement of this $1$-dimensional space inside $\mathfrak{sp}(2)$, the $9$-dimensional normal space is
        \[
        \nu_e(U \cdot e) = \left\{ \begin{pmatrix}  \alpha j & x \\ -\thickbar{x} & v \end{pmatrix} : \alpha \in \C, \ x \in \mathbb{H}, \ v \in \mathrm{Im}(\mathbb{H}) \right\}.
        \]
        The corresponding slice representation is such that the normal space decomposes into invariant subspaces as $\nu_e(U \cdot e) \cong \C \oplus \C^2 \oplus \R^3$, with fixed subspace $\C$, where $U_e^0$ acts on $\C^2$ via the $\SU(2)$ standard representation and it acts on $\R^3$ via the $\SO(3)$ adjoint representation. The restricted action of $U_e^0 = \Delta \Sp(1)_{\low}$ on the $7$-dimensional submodule $\R^3 \oplus \C^2$ is not polar since its cohomogeneity on $\R^3 \oplus \C^2$ is $4$, whereas the sum of the cohomogeneities on $\R^3$ and $\C^2$ is $1 + 1 = 2 \neq 4$. Therefore, by Theorem~\ref{Thm:Dadok}\ref{Thm:Dadok(i)}, the action is not infinitesimally polar.

        \item $\mathfrak{n}_1 = \mathfrak{u}(1), \mathfrak{h}_2 = \mathfrak{su}(2)_\Delta, \mathfrak{n}_2 = \{0\}$. Owing to Goursat's Lemma, the subalgebra $\mathfrak{u}$ arises as
        $$\mathfrak{u} = \{ (\diag(i\alpha, q), \diag(q, q)) : \alpha \in \R, \, q \in \mathfrak{sp}(1) \}.$$ 
        The isotropy subalgebra at $e \in \Sp(2)$ consists of elements where the left and right factors coincide, meaning $\diag(i\alpha, q) = \diag(q, q)$. This equality forces $q = i\alpha$, yielding
        \[
        \mathfrak{u}_e = \{ (\diag(i\alpha, i\alpha), \diag(i\alpha, i\alpha)) : \alpha \in \R \}.
        \]
        Hence, the connected component of the isotropy subgroup is the $1$-dimensional torus $U_e^0 = \Delta T^1_{1,1}$. Moreover, $T_e(U \cdot e) = \{ \diag(i\alpha - q, 0) : \alpha \in \R, \, q \in \mathfrak{sp}(1) \} = \mathfrak{sp}(1)_{\mathrm{up}}$, so the $7$-dimensional normal space is
        \[
        \nu_e(U \cdot e) = \left\{X =  \begin{pmatrix} 0 & -\thickbar{x} \\ x & v \end{pmatrix} : v \in \mathrm{Im}(\mathbb{H}), \, x \in \mathbb{H} \right\}.
        \]
        Recalling that the slice representation at $e$ is the action of conjugation $U_e \curvearrowright \nu_e(U \cdot e)$ and writing $v = ia + z j \in \mathrm{Im}(\mathbb{H})$ (with $a \in \R, z \in \C$) and $x = x_1 + x_2 j \in \mathbb{H}$ (with $x_1, x_2 \in \C$), the conjugation of the blocks of $X \in \nu_e(U \cdot e)$ by $e^{i\theta}$ yields
        \[
        e^{i\theta} v e^{-i\theta} = ia + e^{2i\theta}z j \qandq e^{i\theta} x e^{-i\theta} = x_1 + e^{2i\theta} x_2 j.
        \]
        Hence, $\nu_e(U \cdot e)$ decomposes into $U_e^0$-modules as $\nu_e(U \cdot e) \cong \R^3 \oplus \C_2 \oplus \C_2$. Since $U_e^0$ is a one-dimensional torus, Corollary~\ref{corolario:toro} implies that the restricted representation on $\C_2 \oplus \C_2$ is not polar. Therefore, by Theorem~\ref{Thm:Dadok}\ref{Thm:Dadok(i)}, the slice representation at $e$ is not polar.

        \item $\mathfrak{n}_1 = \mathfrak{u}(1), \mathfrak{h}_2 = \mathfrak{su}(2)_{\low}, \mathfrak{n}_2 = \{0\}$. Owing to Goursat's Lemma, the Lie subalgebra $\mathfrak{u}$ is given by $\mathfrak{u} = \{ (\diag(i\alpha, q), \diag(0, q)) : \alpha \in \R, \, q \in \mathfrak{sp}(1)\}.$
        The isotropy subalgebra at $e \in \Sp(2)$ is given by elements satisfying $\diag(i\alpha, q) = \diag(0, q)$, which forces $\alpha = 0$ while $q \in \mathfrak{sp}(1)$ remains free. Hence, 
        \[\mathfrak{u}_e = \{ (\diag(0, q), \diag(0, q)) : q \in \mathfrak{sp}(1) \},\]
        and the connected component of the isotropy subgroup is precisely $U_e^0 = \Delta\Sp(1)_{\low}$. Since $T_e(U \cdot e) = \{ \diag(i\alpha, 0) : \alpha \in \R \}$, the $9$-dimensional normal space is given by 
        \[
        \nu_e(U \cdot e) = \left\{ X =  \begin{pmatrix} \alpha j & x \\ -\thickbar{x} & v \end{pmatrix} : \alpha \in \C , \ x \in \mathbb{H}, \ v \in \mathrm{Im}(\mathbb{H}) \right\}.
        \]
        The slice representation is given by the conjugation action of $U_e^0$ on $\nu_e(U \cdot e)$. For $h = \diag(1, w) \in \Sp(1)_{\low}$, the conjugation $h X h^{-1}$ yields
        \[
        \begin{pmatrix} 1 & 0 \\ 0 & w \end{pmatrix} \begin{pmatrix} \alpha j & x \\ -\thickbar{x} & v \end{pmatrix} \begin{pmatrix} 1 & 0 \\ 0 & \thickbar{w} \end{pmatrix} = \begin{pmatrix} \alpha j & x\thickbar{w} \\ -w\thickbar{x} & w v \thickbar{w} \end{pmatrix}.
        \]
         It is straightforward to check that $\nu_e(U \cdot e)$ decomposes into $U_e^0$-invariant subspaces as $\nu_e(U \cdot e) \cong \R^2 \oplus \R^3 \oplus \C^2$, where $\R^2$ is the fixed component, $\R^3$ is the submodule where $U_e^0$ acts via the adjoint representation of $\SO(3)$ and $\C^2$ is the subspace on which $U_e^0$ acts via the standard $\SU(2)$-representation. The restricted action of $U_e^0$ on the submodule $\R^3 \oplus \C^2$ is not polar since the sum of the cohomogeneities on $\R^3$ and $\C^2$ is $2$ whereas the cohomogeneity on $\R^3 \oplus \C^2$ is $4$. Therefore, by Theorem~\ref{Thm:Dadok}\ref{Thm:Dadok(i)}, the slice representation at $e$ is not polar.

        \item $\mathfrak{n}_1 = \mathfrak{su}(2)_{\low}, \mathfrak{h}_2 = \mathfrak{t}^2, \mathfrak{n}_2 = \mathfrak{t}^1_{m,n}$. First, note that $\dim \mathfrak{u} = \dim \mathfrak{n}_1 + \dim \mathfrak{n}_2 + 1 = 5$. Since $\mathfrak{h}_2 = \mathfrak{t}^2 \subset \mathfrak{h}_1$ and $\dim \mathfrak{h}_1 = 4$, we have $T_e(U \cdot e) \subset \mathfrak{h}_1$ and hence $\dim \mathfrak{u}_e = \dim \mathfrak{u} - \dim T_e(U \cdot e) \geq 1$. Because $\mathfrak{u}_e \subset \Delta(\mathfrak{h}_1 \cap \mathfrak{h}_2) = \Delta \mathfrak{t}^2$, the connected isotropy subgroup $U_e^0$ is a torus $T^r \subset \Delta T^2$ with $r \in \{1,2\}$. Moreover, $\nu_e(U \cdot e)$ contains the $6$-dimensional $U_e^0$-invariant subspace
        \[
        W \coloneqq \mathfrak{h}_1^\perp = \left\{ \begin{pmatrix} \alpha j & x \\ -\thickbar{x} & 0 \end{pmatrix} : \alpha \in \C, \, x \in \mathbb{H} \right\}.
        \]
        Writing $x = x_1 + x_2 j$ with $x_1, x_2 \in \C$, conjugation by $\diag(e^{i\theta_1}, e^{i\theta_2}) \in T^2$ sends $\alpha j \mapsto e^{2i\theta_1} \alpha j$, $x_1 \mapsto e^{i(\theta_1 - \theta_2)} x_1$, and $x_2 j \mapsto e^{i(\theta_1 + \theta_2)} x_2 j$, so $W \cong \C_{(2,0)} \oplus \C_{(1,-1)} \oplus \C_{(1,1)}$ as a $T^2$-module. If $r = 2$, then $U_e^0 = \Delta T^2$ acts without non-trivial fixed vectors on $W \cong \C^3$ ($r = 2 \neq 3$). If $r = 1$, then $U_e^0 = \Delta T^1_{a,b}$ for some non-zero $(a,b) \in \Z^2$, whose weights on $W$ are $2a$, $a - b$, and $a + b$; at most one of them vanishes, so $U_e^0$ acts without non-trivial fixed vectors on an invariant subspace $\C^s \subset W$ with $s \in \{2,3\}$ ($r = 1 \neq s$). In both cases, up to a finite ineffective kernel, Corollary~\ref{corolario:toro} and Theorem~\ref{Thm:Dadok}\ref{Thm:Dadok(i)} imply that the slice representation at $e$ is not polar, and hence the $U$-action on $\Sp(2)$ is not infinitesimally polar.
        
        \item \label{Subsec455itemi} Consider $\mathfrak{n}_1 = \mathfrak{n}_2 =\mathfrak{su}(2)_{\low}$ and $\mathfrak{h}_2=\mathfrak{u}(1) \oplus \mathfrak{su}(2)$. This choice of Lie subalgebras and its ideals gives rise to a family of isomorphisms $\theta_{m,k}: \mathfrak{h}_1/\mathfrak{n}_1 \to \mathfrak{h}_2/\mathfrak{n}_2$, where $m,k \in \mathbb{Z}$ are non-zero and coprime, inducing a family of subgroups
        $$U_{m,k} = \{ (\diag(z^m, p), \diag(z^k,q)): z \in \Un(1), \  p,q \in \Sp(1) \}.$$
        Noticing that the action of the subgroup $N_1 \subset U_{m,k}$ given by 
        $$N_1 = \{ (\diag(1,p), e): p \in \Sp(1) \} \cong \Sp(1)_{\low} \times \{ e\},$$
        on $\Sp(2)$ has orbit space $\Sp(2)/N_1 \cong S^7$, it turns out that the $U_{m,k}$-action on $\Sp(2)$ has the same orbit space as $U_{m,k}/N_1 \cong \Sp(1) \times \Un(1)$ acting on $S^7$. We claim that $U_{m,k}$ acts infinitesimally polarly on $\Sp(2)$ if and only if $m = \pm k$. In fact, supposing that $m \neq \pm k$, it is easy to check that the slice representation at the point 
        $$g = \begin{pmatrix}
        0 & 1 \\ 1 & 0    
        \end{pmatrix} \in \Sp(2)$$
        corresponds to $(\Un(1), \C^2)$ with weights $m+k$ and $m-k$. Since $m \neq \pm k$, Corollary \ref{corolario:toro} implies that this representation is not polar. On the other hand, if $m = \pm k$, the action of $U_{m,k}/N_1$ on $S^7$ corresponds to the restriction of the polar representation of $U_{m,k}/N_1$ on $\mathbb{H}^{2} \cong \C_0 \oplus \C_{\lambda} \oplus \C^2$, where $\Un(1)$ acts on $\C_\lambda$ with weight $\lambda \neq 0$, $\Sp(1) \times \Un(1)$ acts on $\C^2$ orbit equivalently to the standard $\Un(2)$-representation, and $\C_0$ denotes the fixed component. Consequently, the $U_{m,k}$-action on $\Sp(2)$ is infinitesimally polar and the biquotient $\Sp(2) \biq U_{m,k}$ is a $3$-dimensional Riemannian orbifold.

        \item $\mathfrak{n}_1=\mathfrak{u}(1), \mathfrak{h}_2= \mathfrak{u}(1) \oplus \mathfrak{su}(2)$ and $\mathfrak{n}_2 = \mathfrak{u}(1)$. This quintuple gives rise to the group $U$ given by $U = \{ (\diag(z,p), \diag(w,p)): z,w \in \Un(1), p \in \Sp(1) \},$
    whose action on $\Sp(2)$ has 
    $$U_e = \{ (\diag(z,p), \diag(z,p)): z \in \Un(1), p \in \Sp(1) \} \cong \Un(1) \times \Sp(1)$$
    as the isotropy subgroup at $e$. Notice that the tangent space $T_e (U \cdot e)$ can be identified with $\mathfrak{u}(1)_{\mathrm{up}} \coloneqq \{ \diag(it, 0) : t \in \R \}$ and hence the normal space $\nu_e(U \cdot e )$ is given by
    $$\nu_e (U \cdot e ) = \left\{ \begin{pmatrix} wj & -\thickbar{x} \\ x & v
    \end{pmatrix}: w \in \C, \ v \in \mathrm{Im}(\mathbb{H}), \ x \in \mathbb{H} \right\}.$$
    With some computations, it is easy to check that $\nu_e(U \cdot e)$ decomposes into $U_e$-irreducible subspaces as $\nu_e(U \cdot e) \cong \C \oplus \R^3 \oplus \mathbb{H}$, where $\Un(1)$ acts on $\C$ with weight 2, $\Sp(1)$ acts on $\R^3$ as the $\SO(3)$ $\Ad$-representation, while $\Sp(1)$ and $\Un(1)$ act on $\mathbb{H}$ by left and right multiplication, respectively. Since its cohomogeneity is $5 \neq 1 + 1 + 1$, Theorem~\ref{Thm:Dadok}\ref{Thm:Dadok(i)} implies that the slice representation at $e$ is not polar and $U$ does not act infinitesimally polarly on $\Sp(2)$. 

    \end{enumerate}

  \subsubsection{$\mathfrak{h}_1=\mathfrak{t}^2$}\label{subsec:t2_goursat}

  \begin{enumerate}[label=\textnormal{(\roman*)}, itemsep=0.2cm]
      \item \label{Subsec:456itemi} $\mathfrak{n}_1=\{0\}, \mathfrak{h}_2=\mathfrak{t}^2$ and $\mathfrak{n}_2=\{0\}$. In this case, the isomorphism $\theta : \mathbb{R}^2 \to \mathbb{R}^2$ is uniquely represented by an invertible matrix $A = \left(\begin{smallmatrix} a & b \\ c & d \end{smallmatrix}\right) \in \GL(2,\Q)$. By Goursat's Lemma, the subalgebra $\mathfrak{u}$ has dimension $2$ and consists of pairs $(X,Y)$ where $X = \diag(i \lambda_1, i \lambda_2, -i \lambda_1, -i \lambda_2)$ and 
      \[Y = \diag\left(i(a\lambda_1 + b\lambda_2), i(c\lambda_1 + d\lambda_2), -i(a\lambda_1 + b\lambda_2), -i(c\lambda_1 + d\lambda_2)\right),\]
      generating the $2$-dimensional subgroup $\Delta^\theta T^2 \subset T^2 \times T^2$. Since $X$ and $Y$ both belong to the maximal torus $\mathfrak{t}^2$, they are conjugate in $\Sp(2)$ if and only if they lie in the same orbit under the action of the Weyl group $W(\Sp(2)) \cong D_4$, the dihedral group of order $8$. Let $v_X =(\lambda_1, \lambda_2) \in \R^2$ be the parameter vector generating $X$ and $v_Y=Av_X$. The conjugacy condition yields $v_X = s v_Y$ for some $s \in D_4$, which is equivalent to the homogeneous linear system $(\id-s A)v_X=0$. Hence, if $\det(\id - sA) \neq 0$ for all $s \in D_4$, this forces $v_X=0$, which implies $X=Y=0$. Consequently, the $\Delta^\theta T^2$-action on $\Sp(2)$ is almost free and the biquotient $\Sp(2) \biq \Delta^\theta T^2$ is an $8$-dimensional Riemannian orbifold.

      Otherwise, if $\det(\id - sA) = 0$ for some $s \in D_4$, up to conjugation by the Weyl group we may assume $\det(\id - A) = 0$, so that the connected isotropy subgroup at the identity $e \in \Sp(2)$ is $(\Delta^\theta T^2)_e^0 = \Delta H \cong H$, where $H \subset T^2$ is a non-trivial torus. Since $T_e(\Delta^\theta T^2 \cdot e) \subset \mathfrak{t}^2 \subset \mathfrak{u}(2)$, the normal space $\nu_e(\Delta^\theta T^2 \cdot e)$ contains $\mathfrak{u}(2)^\perp = \mathrm{Sym}(2,\C)$, and by Lemma \ref{lemma:non-polar-u(2)} the slice representation at $e$ is not polar.

      \item $\mathfrak{n}_1 = \mathfrak{t}^1_{k_1,k_2}, \mathfrak{h}_2 = \mathfrak{t}^2, \mathfrak{n}_2 = \mathfrak{t}^1_{m,n}$. First, note that $\dim \mathfrak{u}=\dim \mathfrak{n}_1 + \dim \mathfrak{n}_2+1=3$. Since the tangent space to the orbit at the identity $e \in \Sp(2)$ is a subspace of $\mathfrak{t}^2$, its dimension is at most $2$. Evaluating the dimension of the isotropy subalgebra yields 
      \[
       \dim \mathfrak{u}_e = \dim \mathfrak{u} - \dim T_e(U \cdot e) \geq 3 - 2 = 1.
      \]
      Hence, the connected component of the isotropy subgroup $U_e^0$ is isomorphic to a non-trivial, closed, connected subgroup of $T^2 \subsetneq \Un(2)$. Moreover, since $T_e(U \cdot e) \subset \mathfrak{t}^2 \subset \mathfrak{u}(2)$, its orthogonal complement inside $\mathfrak{sp}(2)$ must contain the entire symmetric block. Consequently, $\Sym(2,\C)$ is an invariant subspace of the normal space $\nu_e(U \cdot e)$. By Lemma~\ref{lemma:non-polar-u(2)}, the restricted action of $U_e^0$ on the invariant subspace $\Sym(2,\C)$ is not polar. Therefore, by Theorem~\ref{Thm:Dadok}\ref{Thm:Dadok(i)}, the $U_e^0$-action on the normal space $\nu_e(U \cdot e)$ is not polar. 

      \item \label{Subsec456itemii}$\mathfrak{n}_1=\mathfrak{t}^1_{k_1,k_2}, \mathfrak{h}_2= \mathfrak{t}^1_{m,n}$ and $\mathfrak{n}_2 = \{0\}$. Since $\gcd(k_1,k_2) =1$, the quotient for the first factor is $\mathfrak{h}_1/\mathfrak{n}_1 \cong \R$, and the corresponding projection $\pi_1 : \mathfrak{t}^2 \to \R$ associates each $X = \diag(ix, iy, -ix, -iy)$ with the real number $\pi_1(X) = -k_2x+k_1y$. For the second factor, $\mathfrak{h}_2/\mathfrak{n}_2 \cong \R$ and the projection $\pi_2 : \mathfrak{t}^1_{m,n} \to \R$ associates each $Y = \diag(imz, inz, -imz, -inz)$ with $\pi_2(Y) = z \in \R$. The isomorphism $\theta$ corresponds to multiplication by a non-zero rational constant $a \in \Q^\ast$, so the fiber product condition $\theta(\pi_1(X)) = \pi_2(Y)$ translates to $z = a(-k_2x+k_1y)$. Hence, this quintuple gives rise to the $2$-dimensional subgroup $U = T^2_a(k_1,k_2,m,n) \subset T^2 \times T^1_{m,n}$ given by
      \[
      T^2_a(k_1,k_2,m,n) = \left\{ \big( \diag(e^{ix}, e^{iy}, e^{-ix}, e^{-iy}), \, \diag(e^{imz}, e^{inz}, e^{-imz}, e^{-inz}) \big) : x,y \in \R \right\},
      \]
where $z = a(-k_2x+k_1y)$. At the identity $e \in \Sp(2)$, the isotropy condition $X=Y$ yields $x = mz$ and $y = nz$. Substituting these into $z = a(-k_2x+k_1y)$ gives $z=a(k_1 n - k_2m)z$, which admits a non-zero solution for $z$ if and only if $k_1 n - k_2 m \neq 0$ and $a = (k_1 n - k_2 m)^{-1}$. More generally, since the Weyl group $W(\Sp(2)) \cong D_4$ acts on $\mathfrak{t}^2$ via signed permutations of the diagonal entries, by Lemma \ref{lemma:almostfree} a non-zero element $(X,Y) \in \mathrm{Lie}(U)$ satisfies $X = \Ad_g Y$ for some $g \in \Sp(2)$ if and only if $a \langle (-k_2, k_1), s(m,n) \rangle = 1$ for some $s \in D_4$. Consequently, if $a \langle (-k_2, k_1), s(m,n) \rangle \neq 1$ for all $s \in D_4$, the $U$-action on $\Sp(2)$ is almost free and the biquotient $\Sp(2) \biq T^2_a(k_1,k_2,m,n)$ is an $8$-dimensional Riemannian orbifold.

Otherwise, up to conjugation by the Weyl group (which permutes and changes the signs of $m$ and $n$), we may assume $k_1 n - k_2 m \neq 0$ and $a = (k_1 n - k_2 m)^{-1}$, so it suffices to analyze the isotropy at the identity $e \in \Sp(2)$. In this case, $x = mz$ and $y = nz$, and the connected isotropy subgroup at $e$ is $U_e^0 = \Delta T^1_{m,n} \cong T^1_{m,n}$. Since $T_e(U \cdot e) = \mathfrak{n}_1 \subset \mathfrak{t}^2 \subset \mathfrak{u}(2)$, the normal space $\nu_e(U \cdot e)$ contains $\mathfrak{u}(2)^\perp = \mathrm{Sym}(2,\C)$, and by Lemma \ref{lemma:non-polar-u(2)} the slice representation at $e$ is not polar.
  
  \end{enumerate}  
  
  \subsubsection{$\mathfrak{h}_1=\mathfrak{t}^1_{m,n}, \mathfrak{n}_1=\mathfrak{n}_2=\{0\}, \mathfrak{h}_2=\mathfrak{t}^1_{r,s}$}\label{subsec:torus_asymm} For this choice of $\mathfrak{h}_1, \mathfrak{h}_2, \mathfrak{n}_1$ and $\mathfrak{n}_2$, the isomorphism $\theta$ is a multiplication by a rational parameter $p/q \in \Q^{\ast}$, such that $\gcd(p,q) = 1$. First, consider $\vert p \vert = \vert q \vert$ and assume that $\vert p \vert = 1$. Up to conjugation, the subgroup $U$ is given by 
  $$U_{m,n}^{r,s} = \{ (\diag(z^{m}, z^{n}, z^{-m}, z^{-n}), \diag(z^{r}, z^{s}, z^{-r}, z^{-s})): z \in \Un(1) \},$$
  which acts almost freely on $\Sp(2)$ if and only if $\{ \vert m \vert , \vert n \vert \} \neq \{ \vert r \vert, \vert s \vert \}$, and the corresponding biquotient $\Sp(2) \biq U_{m,n}^{r,s}$ is a $9$-dimensional Riemannian orbifold. On the other hand, if $\{ \vert m \vert , \vert n \vert \} = \{ \vert r \vert, \vert s \vert \}$, the $U$-action on $\Sp(2)$ corresponds to the diagonal action $\Delta T^{1}_{m,n} \curvearrowright \Sp(2)$, which is not infinitesimally polar, see Subsection~\ref{seccao:diagonal}. Now, assume that $\vert p \vert \neq \vert q \vert$ and note that, in this case, the corresponding subgroup is 
  $$H_{p,q}(m,n,r,s) = \{ (\diag(z^{pm}, z^{pn}, z^{-pm}, z^{-pn}), \diag(z^{qr}, z^{qs}, z^{-qr}, z^{-qs})): z \in \Un(1) \}.$$
  Since $\gcd(m,n) = \gcd(r,s) = 1$, we have $\gcd(|pm|, |pn|) = |p| \neq |q| = \gcd(|qr|, |qs|)$, and hence $\{ \vert pm \vert , \vert pn \vert \} \neq \{ \vert qr \vert, \vert qs \vert \}$, implying that the action of $H_{p,q}(m,n,r,s)$ on $\Sp(2)$ is almost free. Therefore, for all $m,n,r,s \in \Z$ such that $\gcd(m,n) = \gcd(r,s) = 1$, the biquotient $\Sp(2) \biq H_{p,q}(m,n,r,s)$ is a $9$-dimensional Riemannian orbifold, whenever $\vert p \vert \neq \vert q \vert$.

\subsection{Infinitesimally polar actions on $\SO(5)$}\label{subsec:classification_so(5)}

As a direct consequence of Theorem~\ref{thm:main_classification}, we now deduce the classification of all infinitesimally polar isometric actions of compact connected Lie groups on $\SO(5)$ equipped with a bi-invariant metric. 

\medskip
\noindent
\textit{Proof of Corollary \textnormal{\ref{cor:SO5_classification}}.}
Since $\ker \Psi$ is finite and central, the map $U \mapsto \Psi(U)$ defines a bijection (with inverse $\tilde{U} \mapsto (\Psi^{-1}(\tilde{U}))^0$) between closed connected subgroups of $\Sp(2) \times \Sp(2)$ and of $\SO(5) \times \SO(5)$, preserving conjugacy, factor permutation, and non-transitivity. For any $g \in \Sp(2)$ and $\tilde{g} = \psi(g)$, the local isometry $\psi$ is equivariant with respect to $\Psi|_U \colon U \to \Psi(U)$, so $\d\psi_g$ conjugates the slice representation of $U_g^0$ on $\nu_g(U \cdot g)$ to that of $(\Psi(U)_{\tilde{g}})^0$ on $\nu_{\tilde{g}}(\Psi(U) \cdot \tilde{g})$. Because polarity of a linear representation depends only on the identity component of the acting group, $\Psi(U) \curvearrowright \SO(5)$ is infinitesimally polar if and only if $U \curvearrowright \Sp(2)$ is infinitesimally polar. \hfill $\square$

\section{New examples of positively curved Riemannian orbifolds}\label{sec:pos_curv}

In this section, we establish Corollaries~\ref{cor:positive_curvature} and~\ref{cor:positive_curvature_4d}, showing that the almost-free actions of $T^1_{m,n}$ and $\SU(2)_{\Delta}$ on the positively curved Berger space $B^7 = \Sp(2)/\SU(2)_{\irr}$ produce new positively curved Riemannian orbifolds with non-trivial singular locus.

\medskip

\noindent
\textit{Proof of Corollary \textnormal{\ref{cor:positive_curvature}}.} It was shown in Subsection~\ref{section:berger-space}\ref{subsec:berger-o6} that, for any pair of coprime integers $(m,n)$ such that $\{|m|, |n|\} \neq \{1,3\}$, the one-dimensional torus $T^{1}_{m,n}$ acts almost freely on the Berger space $B^7 = \Sp(2)/K$, where $K = \SU(2)_{\irr}$. Consequently, the biquotients $\mathcal{O}^{6}_{m,n} \coloneqq \Sp(2)\biq (T^{1}_{m,n} \times \SU(2)_{\irr})$ form a family of $6$-dimensional Riemannian orbifolds. Furthermore, by O'Neill's formula for Riemannian submersions (see \cite{oneill}), since $B^7$ is positively curved (see \cite{berger}), the Riemannian orbifolds $\mathcal{O}^{6}_{m,n}$ likewise inherit metrics with positive sectional curvature. 

By \eqref{Def:SU(2)irr}, the diagonal matrices in $\SU(2)_{\irr}$ are precisely those of its maximal torus $T^{1}_{3,-1} = \{ \diag(w^3, w^{-1}, w^{-3}, w) : w \in \Un(1) \}$, so the isotropy subgroup of $T^{1}_{m,n}$ at the identity coset $eK \in B^7$ is given by $T^{1}_{m,n} \cap K = T^{1}_{m,n} \cap T^{1}_{3,-1}$. Since $\gcd(m,n) = 1$, an element $Z = \diag(z^m, z^n, z^{-m}, z^{-n}) \in T^{1}_{m,n}$, with $z \in \Un(1)$, belongs to $T^{1}_{3,-1}$ if and only if $w = z^{-n}$ and $z^{m} = w^3 = z^{-3n}$, that is, $z^{m+3n} = 1$. Since $\{|m|, |n|\} \neq \{1,3\}$, we have $m + 3n \neq 0$, and hence the isotropy subgroup at $eK$ is
\[
T^{1}_{m,n} \cap T^{1}_{3,-1} \cong \Z_{\vert m+3n\vert}.
\]
More generally, at the cosets $gK$ with $g \in N_{\Sp(2)}(T^2)$, conjugating $T^1_{m,n}$ by the Weyl group $W(\Sp(2)) \cong D_4$ yields isotropy subgroups isomorphic to $\Z_{\vert m \pm 3n\vert}$ and $\Z_{\vert 3m \pm n\vert}$. To identify the orbifold isotropy groups, note that the ineffective kernel of the $\Sp(2)$-action on $B^7$ is $Z(\Sp(2)) = \{\pm I_4\} \cong \Z_2$. If $m$ and $n$ have opposite parity, then $-I_4 \notin T^{1}_{m,n}$, so $T^{1}_{m,n}$ acts effectively on $B^7$ and $\max(|m|+3|n|, 3|m|+|n|) \geq 3$, giving a non-trivial orbifold isotropy group of order at least $3$. If $m$ and $n$ are both odd, then the ineffective kernel of $T^{1}_{m,n}$ on $B^7$ is $\{\pm I_4\} \cong \Z_2$, and $|m| + 3|n| \geq 4$, so at a suitable coset $gK$ with $g \in N_{\Sp(2)}(T^2)$, the effective isotropy group $\Z_{(|m|+3|n|)/2}$ has order at least $2$. Therefore, in all cases, $\mathcal{O}^{6}_{m,n}$ has a non-trivial singular locus.
\hfill $\square$
\medskip

We now turn to the $4$-dimensional orbifold $\mathcal{Q}^4 = \Sp(2) \biq (\SU(2)_{\Delta} \times \SU(2)_{\irr})$, for which passing to $\SO(5)$ provides a geometric realization as an $\SO(3)$-quotient of a Stiefel manifold and allows us to determine its local isotropy groups, its singular locus, and its underlying topology.

\medskip

\noindent
\textit{Proof of Corollary \textnormal{\ref{cor:positive_curvature_4d}}.} As established in Subsection~\ref{section:berger-space}\ref{subsec:berger-q4}, $H \coloneqq \SU(2)_{\Delta}$ acts almost freely on the Berger space $B^7 = \Sp(2)/K$, where $K \coloneqq \SU(2)_{\irr}$. Thus, $\mathcal{Q}^4 \coloneqq \Sp(2) \biq (H \times K)$ is a $4$-dimensional Riemannian orbifold, which inherits a metric of positive sectional curvature from $B^7$ \cite{berger} via O'Neill's formula \cite{oneill}.

Consider the universal double covering $\psi \colon \Sp(2) \to \SO(5)$ with kernel $\{\pm I_4\}$. Since $-I_4 \in H \cap K$, the element $(-I_4, I_4) \in H \times K$ acts on $\Sp(2)$ as the deck transformation $g \mapsto -g$. Swapping the two factors via the isometry $g \mapsto g^{-1}$, we obtain an orbifold isometry
\[
\mathcal{Q}^4 \cong \SO(3)_{\irr} \backslash \SO(5) / \SO(3)_{\mathrm{up}},
\]
where $\SO(3)_{\mathrm{up}}$ is the upper-left $3 \times 3$ block embedding in $\SO(5)$ and $\SO(3)_{\irr}$ acts irreducibly on $\R^5$. Identifying $\R^5$ with the space $\Sym_0(3,\R)$ of traceless symmetric $3 \times 3$ real matrices endowed with the inner product $\langle X, Y \rangle = \tr(XY)$, the irreducible $\SO(3)_{\irr}$-action corresponds to conjugation $X \mapsto gXg^t$. Moreover, since $\SO(5)$ acts transitively on the Stiefel manifold $V_2(\R^5)$ of orthonormal pairs with isotropy $\SO(3)_{\mathrm{up}}$ at the last two standard basis vectors, we have an $\SO(5)$-equivariant identification $\SO(5)/\SO(3)_{\mathrm{up}} \cong V_2(\Sym_0(3,\R))$. Hence, $\mathcal{Q}^4$ is orbifold-diffeomorphic to $\SO(3) \backslash V_2(\Sym_0(3,\R))$, with $\SO(3)$ acting diagonally by conjugation, $g \cdot (A,B) = (gAg^t, gBg^t)$. As $\SO(3)$ is simple and acts non-trivially, this action is effective, so the local orbifold isotropy group at $[(A,B)] \in \mathcal{Q}^4$ coincides with the isotropy subgroup $\SO(3)_{(A,B)}$.

Let $g \in \SO(3)_{(A,B)} \setminus \{I_3\}$, viewed as a rotation by $\varphi \in (0,\pi]$ about a line $\R v \subset \R^3$ with $v \neq 0$. Since $g$ commutes with $A$ and $B$, the vector $v$ is a common eigenvector of $A$ and $B$. If $\varphi \in (0,\pi)$, any symmetric operator on $v^\perp$ commuting with $g|_{v^\perp}$ is scalar. As $A$ and $B$ are traceless, this would force them to be linearly dependent, a contradiction. Hence $\varphi = \pi$, so $g^2 = I_3$. It follows that $\SO(3)_{(A,B)}$ is an Abelian group of symmetric orthogonal matrices, hence simultaneously diagonalizable in $\SO(3)$ and conjugate to a subgroup of 
$$S(\Z_2^3) = \{\diag(\epsilon_1, \epsilon_2, \epsilon_3) : \epsilon_i = \pm 1, \, \epsilon_1\epsilon_2\epsilon_3 = 1\} \cong \Z_2 \times \Z_2.$$ 
Thus $\SO(3)_{(A,B)}$ is isomorphic to $\{I_3\}$, $\Z_2$, or $\Z_2 \times \Z_2$ and all three possibilities occur. Indeed, the orthonormal pair $A_1 = \frac{1}{\sqrt{6}}\diag(-1,-1,2)$ and $B_1 = \frac{1}{\sqrt{2}}\diag(1,-1,0)$ spans the space of traceless diagonal matrices, so $\SO(3)_{(A_1,B_1)} = S(\Z_2^3) \cong \Z_2 \times \Z_2$. For 
$$A_2 = \frac{1}{\sqrt{2}}\begin{pmatrix} 0 & 1 & 0 \\ 1 & 0 & 0 \\ 0 & 0 & 0 \end{pmatrix}$$
and $B_2 = B_1$, both matrices are fixed by $\diag(-1,-1,1)$ while $[A_2, B_2] \neq 0$, so they are not simultaneously diagonalizable and $\SO(3)_{(A_2,B_2)} = \langle \diag(-1,-1,1) \rangle \cong \Z_2$. On the other hand,
$$A_3 = \frac{1}{\sqrt{2}}\begin{pmatrix} 0 & 0 & 1 \\ 0 & 0 & 0 \\ 1 & 0 & 0 \end{pmatrix} \ \text{and} \ B_3 = \frac{1}{\sqrt{2}}\begin{pmatrix} 0 & 0 & 0 \\ 0 & 0 & 1 \\ 0 & 1 & 0 \end{pmatrix}$$
share no common eigenvector (the eigenvectors of $A_3$ are multiples of $(1,0,\pm 1)^t$ and $(0,1,0)^t$, none of which is an eigenvector of $B_3$), so $\SO(3)_{(A_3,B_3)} = \{I_3\}$.

From the homotopy exact sequence of the fibration $\SU(2)_{\irr} \hookrightarrow \Sp(2) \to B^7$, the space $B^7$ is simply connected, and in particular orientable. Since $\mathcal{Q}^4 = (H/\{\pm I_4\}) \backslash B^7$ with $H/\{\pm I_4\} \cong \SO(3)$ connected and acting effectively and almost freely, $\mathcal{Q}^4$ is a closed orientable Riemannian orbifold whose orbifold fundamental group fits into $\pi_1(B^7) \to \pi_1^{\mathrm{orb}}(\mathcal{Q}^4) \to \pi_0(\SO(3))$, whence $\pi_1^{\mathrm{orb}}(\mathcal{Q}^4) = \{0\}$ (see \cite{caramello}). Since any good orbifold is the quotient of its smooth universal covering manifold by its orbifold fundamental group, $\Sigma(\mathcal{Q}^4) \neq \emptyset$ implies that $\mathcal{Q}^4$ is a bad Riemannian orbifold.

Now consider the circle $S^1 \coloneqq \{\diag(I_3, R(\theta)) : \theta \in \R\} \subset \SO(5)$ centralizing $\SO(3)_{\mathrm{up}}$, where $R(\theta) \in \SO(2)$ denotes rotation by $\theta$. Because the metric on $\SO(5)$ is bi-invariant, right multiplication by $S^1$ on $\SO(5)/\SO(3)_{\mathrm{up}}$ induces an isometric $S^1$-action on $\mathcal{Q}^4$ rotating $(A,B)$ within $\Span\{(A,B)\}$. A point $[(A,B)] \in \mathcal{Q}^4$ is fixed by $S^1$ if and only if some circle subgroup of $\SO(3)$ preserves $\Span\{(A,B)\}$ and acts non-trivially on it. Up to conjugation in $\SO(3)$, this circle is $\{\diag(R(\varphi), 1) : \varphi \in \R\}$, under which $\Sym_0(3,\R)$ decomposes orthogonally as $\R A_1 \oplus W_2 \oplus W_1$, where $\R A_1$ is fixed, $W_2 = \Span\{(A_2, B_2)\}$ is rotated by $-2\varphi$, and $W_1 = \Span\{(A_3, B_3)\}$ is rotated by $\varphi$. Since $W_2$ and $W_1$ are irreducible $2$-dimensional representations with distinct non-zero weights (while $\R A_1$ is trivial), they are the only two $2$-planes invariant under this circle subgroup. Moreover, the rotations $\diag(R(\varphi), 1)$ together with $\diag(1,-1,-1)$ act as $\On(2)$ on both $W_2$ and $W_1$, so every orthonormal pair spanning $W_2$ (resp.\ $W_1$) lies in the $\SO(3)$-orbit of $(A_2, B_2)$ (resp.\ $(A_3, B_3)$). Since $[(A_2, B_2)]$ has isotropy $\Z_2$ and $[(A_3, B_3)]$ is regular, $\mathrm{Fix}(S^1, \mathcal{Q}^4) = \{[(A_2, B_2)], [(A_3, B_3)]\}$ consists of exactly two isolated points. Passing to the effective quotient $T^1$, Theorem~1.1 of \cite{harveysearle} shows that the underlying topological space of $\mathcal{Q}^4$ is $T^1$-equivariantly homeomorphic either to $S^4$ with a linear circle action or to a weighted complex projective space with a circle action induced by a linear $T^2$-action on $S^5$. However, as shown in \cite[Section~4.4]{harveysearle}, any effective circle action on a weighted complex projective space induced by a linear $T^2$-action on $S^5$ either acts fixed-point-homogeneously (so its fixed-point set has a $2$-dimensional component) or has exactly three isolated fixed points. Hence, the weighted complex projective space case is excluded, and therefore the underlying topological space of $\mathcal{Q}^4$ is homeomorphic to the sphere $S^4$.

Finally, we determine the singular locus $\Sigma(\mathcal{Q}^4)$. Observe that a point $(A,B) \in V_2(\Sym_0(3,\R))$ projects into $\Sigma(\mathcal{Q}^4)$ if and only if it is fixed by some subgroup $\Z_2 \subset \SO(3)$. Up to conjugation, we may assume that this $\Z_2$ is generated by $\diag(-1,-1,1)$, whose fixed subspace in $\Sym_0(3,\R)$ is the orthogonal sum $W_2 \oplus \R A_1$. Accordingly, any orthonormal pair $(A,B)$ fixed by $\diag(-1,-1,1)$ decomposes as
\[
A = A_0 + a A_1 \qandq B = B_0 + b A_1,
\]
where $A_0, B_0 \in W_2$ and $a = \langle A, A_1 \rangle$, $b = \langle B, A_1 \rangle \in \R$. Since $A$ and $B$ are orthonormal, $W_2 \perp A_1$, and $|A_1| = 1$, we have $|A_0|^2 = 1-a^2$, $|B_0|^2 = 1-b^2$, and $\langle A_0, B_0 \rangle = -ab$, the Cauchy--Schwarz inequality for $A_0$ and $B_0$ then implies that $a^2 + b^2 \leq 1$. Furthermore, since $a$ and $b$ determine the inner products of $A_0, B_0 \in W_2$, they determine $(A_0, B_0)$ uniquely up to rotations and reflections on $W_2$ realized by the elements $\diag(R(\varphi), 1)$ and $\diag(1,-1,-1)$ of $\SO(3)$ (which commute with $\diag(-1,-1,1)$ and fix $A_1$). It follows that the complete closed disk $D^2$ given by $a^2 + b^2 \leq 1$ occurs. In particular, equality $a^2 + b^2 = 1$ occurs if and only if $A_0, B_0 \in W_2$ are linearly dependent, equivalently $[A,B] = [A_0, B_0] = 0$, so that $A$ and $B$ are simultaneously diagonalizable and hence fixed by a conjugate of $S(\Z_2^3)$. Therefore, the interior points of $D^2$ have local isotropy group exactly $\langle \diag(-1,-1,1) \rangle \cong \Z_2$, whereas the boundary points have local isotropy group conjugate to $S(\Z_2^3)$. Any $g \in \SO(3)$ mapping one interior point of $D^2$ to another must normalize $\langle \diag(-1,-1,1) \rangle$ and thus lies in $S(\On(2) \times \On(1))$, which fixes $A_1$ and preserves $(a,b) = (\langle A, A_1 \rangle, \langle B, A_1 \rangle)$; hence distinct interior points of $D^2$ represent distinct points of $\Sigma(\mathcal{Q}^4)$. On the boundary circle $a^2 + b^2 = 1$, two pairs lie in the same $\SO(3)$-orbit if and only if they are related by $N_{\SO(3)}(S(\Z_2^3))/S(\Z_2^3) \cong S_3$, acting by permuting the three diagonal entries. The transposition swapping the first two entries fixes $A_1$ (and therefore acts trivially on $(a,b)$), whereas the $3$-cycle permuting the three entries acts on the boundary circle $a^2 + b^2 = 1$ by rotation of angle $2\pi/3$. Consequently, $\Sigma(\mathcal{Q}^4)$ is homeomorphic to the space obtained by attaching the disk $D^2$ to the circle $S^1$ via the map $z \mapsto z^3$. \hfill $\square$

\end{document}